\documentclass{amsart}

\usepackage{amssymb, amsmath, amsthm}
\usepackage{mathrsfs}
\usepackage[fleqn,tbtags]{mathtools}
\usepackage[shortlabels]{enumitem}
\usepackage{tikz}
\usetikzlibrary{arrows.meta}
\usepackage{xcolor}
\usepackage[colorlinks,linkcolor={blue},citecolor={blue},urlcolor={blue!50!black}]{hyperref}

\let\le\leqslant
\let\leq\leqslant
\let\ge\geqslant
\let\geq\geqslant

\theoremstyle{plain}
\newtheorem{theorem}{Theorem}[section]
\theoremstyle{remark}
\newtheorem{remark}[theorem]{Remark}
\newtheorem{definition}{Definition}[section]
\theoremstyle{plain}
\newtheorem{proposition}[theorem]{Proposition}
\newtheorem{corollary}[theorem]{Corollary}
\newtheorem{lemma}[theorem]{Lemma}
\numberwithin{equation}{section}
\newcommand{\D}{\{-1,1\}}
\newcommand{\E}{\mathbb E}
\newcommand{\cR}{\mathcal R}
\newcommand{\Rad}{\operatorname{Rad}}

\newcommand{\Dom}{\mathsf{D}}
\newcommand{\calF}{\mathscr F}
\newcommand{\calG}{\mathscr G}
\newcommand{\calL}{\mathscr L}
\newcommand{\dd}{\,\mathrm d}
\newcommand{\sgn}{\operatorname{sgn}}
\newcommand{\n}{\Vert}

\title[Stein's inequality and dimension-free \(R\)-sectoriality]{The vector-valued Stein inequality and dimension-free \(R\)-sectoriality}

\author{Jan van Neerven}

\address{Delft Institute of Applied Mathematics\\
Delft University of Technology\\P.O. Box 5031\\2600 GA Delft\\The Netherlands}
\email{J.M.A.M.vanNeerven@tudelft.nl}

\begin{document}

\begin{abstract}
For a complex Banach space \(X\), we characterise dimension-free \(R\)-sectoriality of the Walsh operator, the Ornstein--Uhlenbeck operator, the Laplacian, and the Poisson operator by the vector-valued Stein inequality. For each family, the dimension-free \(R\)-bound of the resolvents on the negative real axis equals the Stein constant. Consequently, the resolvent estimate on the negative real axis already implies uniform \(R\)-sectoriality at an angle strictly smaller than \(\pi/2\). We furthermore prove that the full and dyadic Stein properties are equivalent and that they are independent of the exponent \(1<p<\infty\) entering their definition.
\end{abstract}



\maketitle
\enlargethispage{2pt}

\section{Introduction and main results}
\label{sec:introduction}

The extension of scalar \(L^p\)-estimates to Banach space-valued functions often depends on the geometry of the range space; indeed, this is one of the principal themes of the monographs \cite{HNVW1}--\cite{HNVW3}. For the Walsh operators on \(L^p(\D^d;X)\) and the Ornstein--Uhlenbeck operators on \(L^p(\mathbb R^d,\gamma_d;X)\), with \(\gamma_d\) the standard Gaussian measure on \(\mathbb R^d\), Pisier's analyticity theorem \cite{Pisier} gives a common sectorial angle strictly smaller than \(\pi/2\) and dimension-free sectoriality bounds if and only if \(X\) is \(K\)-convex. If \(X\) is UMD, the Walsh and Ornstein--Uhlenbeck operators admit a bounded \(H^\infty\)-calculus of angle strictly smaller than \(\pi/2\) on \(L^p(\mathbb R^d,\gamma_d;X)\), with dimension-free bounds \cite[Theorem~1.1]{Neerven-OU}. Conversely, boundedness of the \(H^\infty\)-calculus for the Ornstein--Uhlenbeck operator on \(L^p(\mathbb R,\gamma;X)\) already implies that \(X\) is UMD; see \cite[Proposition~2.2]{BetancorEtAl} and \cite[Theorem~9.1]{Neerven-OU}.

By a result of Kalton and Weis \cite[Theorem~5.3]{KaltonWeis}, see also \cite[Theorem~10.3.4(2)]{HNVW2}, if a densely defined sectorial operator \(A\) on a Banach space \(E\) with the triangular contraction property has a bounded \(H^\infty\)-calculus, then it is \(R\)-sectorial and \(\omega_R(A)\leq\omega_{H^\infty}(A)\). In turn, \(R\)-sectoriality implies sectoriality with \(\omega(A)\leq\omega_R(A)\).

The combination of these results leads to the natural question whether dimension-free \(R\)-sectoriality of the vector-valued Walsh and Ornstein--Uhlenbeck operators admits a characterisation in terms of a Banach space property of \(X\). The main result of this paper answers this question affirmatively, by showing that the validity of Bourgain's vector-valued Stein inequality for \(X\) is a necessary and sufficient condition for dimension-free \(R\)-sectoriality of both operator families. The same characterisation holds for the Euclidean Laplacian and the Poisson operator. For each of these four families, the dimension-free \(R\)-bound of the resolvents on the negative half-line equals the optimal constant in Stein's inequality for \(X\). To complete this circle of ideas, we show that the vector-valued Stein inequality implies the triangular contraction property.

Using Rota's dilation \cite{Rota}, dimension-free estimates for the \(R\)-bound of the heat semigroup in terms of the \(\mathrm{UMD}^+\) constant of \(X\) were obtained by Hyt\"onen, Li, and Naor \cite[Proposition~20 and Corollary~21]{HytonenLiNaor}. Laplace averaging gives the corresponding bound for the heat resolvents on the negative real axis.
For order-continuous Banach function spaces, Kalton, Lorist and Weis \cite[Theorem~2.4.9]{KaltonLoristWeis} characterised the UMD property by \(\ell^2\)-sectoriality of the Laplacian. In the discussion following that theorem, they obtained a corresponding result for the Poisson generator on Banach function spaces of finite cotype and asked whether \(R\)-sectoriality of the Poisson generator implies UMD on general Banach spaces. Our result show that the answer to this question is affirmative if and only if the validity of the Stein inequality implies the UMD property.

\smallskip
The paper is organised as follows. The notion of $R$-boundedness is introduced in Section \ref{sec:preliminaries}. In Section \ref{sec:Stein}, we investigate the vector-valued Stein inequality as a Banach space property. We show the equivalence with its dyadic counterpart and prove its \(p\)-independence. The main result, which characterises the Stein property in terms of dimension-free \(R\)-sectoriality estimates, is proved in Section \ref{sec:dimensionfree}.

\section{$R$-Boundedness}\label{sec:preliminaries}

We begin by introducing the notions of \(R\)-boundedness and defining the Stein property.

Throughout the paper, we work over the complex scalar field. We denote by \((\varepsilon_j)_{j\ge 1}\) a {\em Steinhaus sequence}, i.e., a sequence of independent variables uniformly distributed on the unit circle in the complex plane. We think of such sequences as the complex-valued analogue of Rademacher sequences, and indeed all results in the paper remain true for Rademacher sequences. In many places, however, this would introduce an additional constant \(2\). This choice of randomisation variables is harmonised with the choice of scalar field, which is complex throughout this paper; see \cite[Section 6.1.a]{HNVW2} for more on this.

The {\em \(R\)-bound} of a family \(\mathscr T\subseteq\calL(E)\), where \(E\) is a Banach space (typically, \(E = L^p(X)\)) is the least \(C\in[0,\infty]\) such that
\[
 \Bigl(\E\Bigl\|\sum_{j=1}^m\varepsilon_jT_jx_j\Bigr\|^2\Bigr)^{1/2}
 \leq C
 \Bigl(\E\Bigl\|\sum_{j=1}^m\varepsilon_jx_j\Bigr\|^2\Bigr)^{1/2}
\]
holds for all finite choices \(T_j\in\mathscr T\) and \(x_j\in E\), with repetitions among the \(T_j\) allowed. This least constant is denoted by \(\cR_E(\mathscr T)\).

\begin{definition} We say that \(\mathscr T\subseteq \calL(E)\) is \emph{\(R\)-bounded} if  \(\cR_E(\mathscr T) < \infty\).
\end{definition}

Every \(R\)-bounded family is uniformly bounded, with \(\sup_{T\in\mathscr T}\n T\n \le \cR_E(\mathscr T)\). For Hilbert spaces \(E\), the notion of \(R\)-boundedness reduces to uniform boundedness. Finite families \(\mathscr T\) are \(R\)-bounded, with \(R\)-bound \(\cR_E(\mathscr T) \le \sum_{T\in \mathscr T}\n T\n\). A comprehensive treatment of \(R\)-boundedness is given in \cite[Chapter 8]{HNVW2}.

We write \(\Rad_m(E)\) for the {\em Rademacher space} consisting of all \(E\)-valued sums \(\sum_{j=1}^m\varepsilon_jx_j\), equipped with the norm
\[
 \Bigl\|\sum_{j=1}^m\varepsilon_jx_j\Bigr\|_{\Rad_m(E)}
 :=\Bigl(\E\Bigl\|\sum_{j=1}^m\varepsilon_jx_j\Bigr\|_E^2\Bigr)^{1/2}.
\]
With this notation, for a family \((T_j)_{j=1}^m\) in \(\calL(E)\), with repetitions allowed, we put \(\mathscr T:=\{T_1,\ldots,T_m\}\) and have
\[
 \cR_E(\mathscr T) = \n T\n_{\calL(\Rad_m(E))},
\]
where \(T(\sum_{j=1}^m\varepsilon_jx_j) = \sum_{j=1}^m\varepsilon_jT_j x_j\).

\section{The Stein property}\label{sec:Stein}

The scalar conditional-expectation inequality underlying the Stein property is due to Stein \cite{Stein-LP}. To formulate its vector-valued version, we recall that if \((S,\calF,\mu)\) is a probability space and \(\calG\subseteq \calF\) is a sub-\(\sigma\)-algebra, then for each \(1\le p\le \infty\) the conditional expectation operator \(\E_{\calG}\) is a contraction on \(L^p(S;X)\) for every Banach space \(X\).
In what follows we fix $1<p<\infty$.

A {\em filtration} is an increasing family of sub-\(\sigma\)-algebras. Loosely speaking, the Stein property asks for the \(R\)-boundedness of the vector-valued conditional expectation operators associated with every filtration. The precise definition is as follows.

\begin{definition}
We say that a Banach space \(X\) has the \emph{\(p\)-Stein property} if there exists a finite constant \(C_{p,X}\) such that for every probability space \((S,\calF,\mu)\) and every finite filtration \(\calF_0\subseteq \dots\subseteq\calF_n\subseteq \calF\) with \(n\geq1\), the associated family of conditional expectation operators \(\{\E_{\calF_0},\dots,\E_{\calF_n}\}\) on \(L^p(S;X)\) satisfies
\[
 \cR_{L^p(S;X)}(\{\E_{\calF_0},\dots,\E_{\calF_n}\}) \le C_{p,X}.
\]
We write \(\mathfrak s_{p,X}\) for the least such constant, with value \(+\infty\) if no finite constant exists.
\end{definition}

A Banach space has the {\em UMD\(_p^+\) property} \cite{Garling}, see also \cite[Section 4.2]{HNVW1}, if
there exists a finite constant \(C_{p,X}\) such that whenever \(S\) is a probability space and
\((d_j)_{j=1}^n\) is a finite martingale difference sequence in \(L^p(S;X)\), then
\[\E\Bigl\|\sum_{j=1}^n \varepsilon_jd_j\Bigr\|_{L^p(S;X)}^p\leq C_{p,X}^p \E \Bigl\|\sum_{j=1}^n d_j\Bigr\|_{L^p(S;X)}.\]
We then write \(\beta_{p,X}^+\) for the least such constant.
By Bourgain's vector-valued Stein inequality, UMD\(_p^+\) spaces (in particular, UMD spaces) have the $p$-Stein property \cite{Bourgain86} and \(\mathfrak s_{p,X}\le \beta_{p,X}^+\); see also \cite[Theorem 4.2.23]{HNVW1}.
By a result of Eskenazis \cite[Lemma 4]{Eskenazis}, every Banach space with the $p$-Stein property is \(K\)-convex.
These results place the $p$-Stein property between \(\mathrm{UMD}^+\) and \(K\)-convexity.

Banach function spaces with the $p$-Stein property are UMD. Indeed, their \(K\)-convexity implies finite cotype, and in a Banach lattice with finite cotype, estimates for Steinhaus sums are equivalent to the corresponding lattice square function estimates. Thus Stein's inequality gives the \(\ell^2\)-boundedness of the dyadic conditional expectations used in \cite[proof of Theorem~2.4.6]{KaltonLoristWeis}; in that proof, UMD\(^+\) is used to obtain precisely this bound. The subsequent weighted argument converts it into uniform martingale transform bounds and hence UMD.

\subsection{Equivalence of the Stein property and the dyadic Stein property}

On the discrete cube \(\D^d\) with its uniform probability measure, the operator \(\mathcal E_n^{(d)}\) denotes conditional expectation onto the first \(n\) coordinates.
For \(1<p<\infty\), the \emph{\(p\)-dyadic Stein property} is defined by restricting the definition of the Stein property to the filtration \(\mathscr F_0^{(d)}\subseteq\cdots\subseteq\mathscr F_d^{(d)}\), where \(\mathscr F_j^{(d)}\) is the \(\sigma\)-algebra generated by the first \(j\) coordinate maps on \(\D^d\). The corresponding least constant is denoted by \(\mathfrak s_{p,X}^{\rm dyad}\), again taking it to be \(+\infty\) if no finite constant exists.

It is clear from the definition that if a Banach space \(X\) has the \(p\)-Stein property, then it has the \(p\)-dyadic Stein property, with bound \( \mathfrak s_{p,X}^{\rm dyad} \le \mathfrak s_{p,X}\). The aim of this subsection is to prove the following converse.

\begin{theorem}[\(p\)-Stein \(=\) \(p\)-dyadic Stein]
\label{thm:full-dyadic-Stein}
Let \(X\) be a Banach space and let \(1<p<\infty\). Then the following assertions are equivalent:
\begin{enumerate}[\rm(1)]
\item \(X\) has the \(p\)-Stein property;
\item \(X\) has the \(p\)-dyadic Stein property.
\end{enumerate}
Moreover, the best constants in the respective definitions are equal:
\[
 \mathfrak s_{p,X}=\mathfrak s^{\rm dyad}_{p,X}.
\]
\end{theorem}

The reduction of general martingale inequalities to dyadic, or Paley--Walsh, models goes back to Maurey; see \cite[Section~3.6]{HNVW1}.

\begin{proof}
Only (2)\(\Rightarrow\)(1) needs proof. We therefore assume that \(C:=\mathfrak s^{\rm dyad}_{p,X}<\infty\).

\smallskip
\emph{Step 1: Reduction to a finite tree.} Let \((S,\mathscr F,\mu)\) be a probability space, let \(\mathscr F_0\subseteq\cdots\subseteq\mathscr F_n\) be a finite filtration, and let \(f_1,\ldots,f_m\in L^p(S;X)\). For each \(\ell\in\{1,\ldots,m\}\), let \(\tau(\ell)\in\{0,\ldots,n\}\); these indices need not be distinct. We write \(\E_j:=\E(\,\cdot\,|\mathscr F_j)\).

The tower property for conditional expectations gives \(\E_{\tau(\ell)}\E_n f_\ell=\E_{\tau(\ell)}f_\ell\). Moreover, every \(\E_n\) acts contractively on their Steinhaus sums, i.e., for all \(f_1,\ldots,f_m\in L^p(S;X)\) we have
\[
 \Bigl\|\sum_{\ell=1}^m\varepsilon_\ell\E_n f_\ell\Bigr\|_{L^p(S;X)}
 \leq
 \Bigl\|\sum_{\ell=1}^m\varepsilon_\ell f_\ell\Bigr\|_{L^p(S;X)}.
\]
Replacing the functions \(f_\ell\) by \(\E_n f_\ell\), we may therefore assume that all \(f_\ell\) are \(\mathscr F_n\)-measurable.

By density, it is enough to prove the Stein estimate for \(\mathscr F_n\)-measurable simple functions \(f_1,\ldots,f_m\). We now fix simple functions \(f_1,\ldots,f_m\) and a number \(\eta>0\). For each \(j<n\), the finitely many functions \(\E_jf_1,\ldots,\E_jf_m\) admit \(\mathscr F_j\)-measurable simple approximations \(v_{j,1},\ldots,v_{j,m}\) satisfying
\[
 \|\E_jf_\ell-v_{j,\ell}\|_{L^p(S;X)}<\eta,
 \qquad 1\leq\ell\leq m.
\]
The sets occurring in their simple-function representations generate a finite \(\sigma\)-algebra \(\mathscr H_j\subseteq\mathscr F_j\) with respect to which all \(v_{j,\ell}\) are measurable. Similarly, there is a finite \(\sigma\)-algebra \(\mathscr H_n\subseteq\mathscr F_n\) with respect to which all \(f_\ell\) are measurable. The \(\sigma\)-algebras
\[
 \mathscr G_j:=\bigvee_{i=0}^j\mathscr H_i,
 \qquad 0\leq j\leq n,
\]
are finite and nested, and satisfy \(\mathscr G_j\subseteq\mathscr F_j\).

Since \(v_{j,\ell}\) is \(\mathscr G_j\)-measurable, the tower property gives, for \(j<n\),
\[
 \E(f_\ell|\mathscr G_j)-\E_jf_\ell
 =
 \E(\E_jf_\ell-v_{j,\ell}\,|\,\mathscr G_j)
 -(\E_jf_\ell-v_{j,\ell}).
\]
Consequently,
\[
 \|\E(f_\ell|\mathscr G_j)-\E_jf_\ell\|_{L^p(S;X)}
 \leq 2\|\E_jf_\ell-v_{j,\ell}\|_{L^p(S;X)}
 <2\eta.
\]
For \(j=n\), both conditional expectations equal \(f_\ell\). It follows that replacing \(\E_{\tau(\ell)}f_\ell\) by \(\E(f_\ell|\mathscr G_{\tau(\ell)})\) changes the random sum by at most \(2m\eta\) in \(\Rad_m(L^p(S;X))\).

Thus, if the Stein estimate holds with constant \(C\) for every finite probability tree, then
\[
 \Bigl\|\sum_{\ell=1}^m\varepsilon_\ell
                 \E_{\tau(\ell)}f_\ell
       \Bigr\|_{\Rad_m(L^p(S;X))}
 \leq
 C\Bigl\|\sum_{\ell=1}^m\varepsilon_\ell f_\ell
       \Bigr\|_{\Rad_m(L^p(S;X))}
 +2m\eta.
\]
Letting \(\eta\downarrow0\) proves the required estimate. It remains to establish the finite-tree estimate.

\smallskip
\begingroup
\def\SteinFigureInManuscript{}
%
%
\ifdefined\SteinFigureInManuscript\else
\documentclass[10pt,a4paper]{amsart}
\usepackage[T1]{fontenc}
\usepackage{lmodern,amsmath,amssymb,mathrsfs}
\usepackage[margin=22mm]{geometry}
\usepackage{tikz}
\usetikzlibrary{arrows.meta}
\pagestyle{empty}
\begin{document}
\fi

\definecolor{steinInk}{HTML}{213443}
\definecolor{steinA}{HTML}{176B78}
\definecolor{steinB}{HTML}{9B581E}
\definecolor{steinRule}{HTML}{C6CDD2}

\begin{figure}[p]
\centering
\begin{minipage}{12.5cm}
\centering
\begin{tikzpicture}[
  x=1cm,y=1cm,
  font=\fontsize{9}{11}\selectfont,
  text=steinInk,
  >={Stealth[length=1.5mm]},
  line/.style={draw=steinInk!65,line width=.55pt,->},
  prob/.style={fill=white,inner sep=2pt,font=\fontsize{9}{10}\selectfont},
  vertex/.style={draw=steinInk!50,fill=white,rounded corners=2pt,
                 minimum width=8mm,minimum height=6mm,inner sep=2pt},
  leaf/.style={vertex,minimum width=13mm,minimum height=8mm,align=center},
  panel/.style={anchor=west,font=\bfseries\fontsize{10}{12}\selectfont},
  note/.style={anchor=west,font=\fontsize{9}{10.5}\selectfont},
  rowval/.style={font=\fontsize{10}{12}\selectfont}
]
\path[use as bounding box] (0,5.5) rectangle (12.5,18.1);

\node[panel] at (0,17.9) {(a)\quad Replace the transition probabilities};
\node[anchor=east] at (12.5,17.9) {$n=1,\quad q=2$};
\node[note] at (0,17.45)
  {On each edge: original probability $\vartheta\ \longmapsto\ \vartheta_2=a_2/4$.};

\node[vertex] (S) at (6.25,16.85) {$S$};
\node[vertex,draw=steinA,fill=steinA!8] (A) at (3.0,15.65) {$A$};
\node[vertex,draw=steinB,fill=steinB!8] (B) at (9.2,15.65) {$B$};
\node[anchor=west,font=\fontsize{8}{10}\selectfont] at (0,16.85) {root: level $-1$};
\node[anchor=west,font=\fontsize{8}{10}\selectfont] at (0,15.65) {$\mathscr G_0$};
\node[anchor=west,font=\fontsize{8}{10}\selectfont] at (0,14.25) {$\mathscr G_1$};

\node[leaf,draw=steinA,fill=steinA!7] (Aone) at (1.55,14.25) {$A_1$\\$f=0$};
\node[leaf,draw=steinA,fill=steinA!18] (Atwo) at (4.05,14.25) {$A_2$\\$f=4$};
\node[leaf,draw=steinB,fill=steinB!6] (Bone) at (7.1,14.25) {$B_1$\\$f=-4$};
\node[leaf,draw=steinB,fill=steinB!15] (Btwo) at (9.2,14.25) {$B_2$\\$f=2$};
\node[leaf,draw=steinB,fill=steinB!25] (Bthree) at (11.3,14.25) {$B_3$\\$f=8$};

\draw[line,draw=steinA] (S) -- node[prob] {$\frac23\mapsto\frac34$} (A);
\draw[line,draw=steinB] (S) -- node[prob] {$\frac13\mapsto\frac14$} (B);
\draw[line,draw=steinA] (A) -- node[prob] {$\frac13\mapsto\frac14$} (Aone);
\draw[line,draw=steinA] (A) -- node[prob] {$\frac23\mapsto\frac34$} (Atwo);
\draw[line,draw=steinB] (B) -- node[prob,pos=.75,yshift=2mm] {$\frac15\mapsto\frac14$} (Bone);
\draw[line,draw=steinB] (B) -- node[prob] {$\frac12\mapsto\frac24$} (Btwo);
\draw[line,draw=steinB] (B) -- node[prob,pos=.75,yshift=2mm] {$\frac3{10}\mapsto\frac14$} (Bthree);

\node[note] at (0,13.45) {$\mu$};
\node[note] at (0,12.94) {$\mu_2$};
\foreach \xx/\old/\new in {1.55/{2/9}/{3/16},4.05/{4/9}/{9/16},7.1/{1/15}/{1/16},9.2/{1/6}/{2/16},11.3/{1/10}/{1/16}} {
  \node at (\xx,13.45) {$\old$};
  \node at (\xx,12.94) {$\new$};
}
\draw[steinRule] (0,12.55) -- (12.5,12.55);

\node[panel] at (0,12.17) {(b)\quad Sample the tree with two blocks of two signs};
\node[note] at (0,11.72)
  {$b_1=(x_1,x_2),\quad b_2=(x_3,x_4),\qquad --,-+,+-,++\ \longleftrightarrow\ 1,2,3,4$.};
\node at (5.15,11.28) {second block $b_2$: choose a child};
\node[align=center,font=\fontsize{9}{10}\selectfont] at (9.35,10.98)
  {current\\vertex};
\node[align=center,font=\fontsize{9}{10}\selectfont] at (11.35,10.98)
  {row average\\$\mathcal E_2^{(4)}J_2f$};

\foreach \col/\sg in {0/{--},1/{-+},2/{+-},3/{++}} {
  \node at ({2.825+1.55*\col},10.86) {$\sg$};
}
\node[rotate=90] at (.3,8.65) {first block $b_1$: choose $A$ or $B$};
\foreach \row/\sg in {0/{--},1/{-+},2/{+-},3/{++}} {
  \node[rowval] at (1.42,{10.15-\row}) {$\sg$};
}

\foreach \row in {0,1,2} {
  \foreach \col/\atom/\val/\shade in {0/A_1/0/7,1/A_2/4/18,2/A_2/4/18,3/A_2/4/18} {
    \fill[steinA!\shade] ({2.05+1.55*\col},{9.65-\row}) rectangle ++(1.55,1);
    \node[align=center] at ({2.825+1.55*\col},{10.15-\row}) {$\atom$\\$f=\val$};
  }
  \draw[line,draw=steinA] (8.43,{10.15-\row}) -- (8.94,{10.15-\row});
  \node[vertex,draw=steinA,fill=steinA!8] at (9.35,{10.15-\row}) {$A$};
  \draw[line,draw=steinA] (9.86,{10.15-\row}) -- (10.61,{10.15-\row});
  \node[rowval,text=steinA] at (11.35,{10.15-\row}) {$3$};
}
\foreach \col/\atom/\val/\shade in {0/B_1/-4/6,1/B_2/2/15,2/B_2/2/15,3/B_3/8/25} {
  \fill[steinB!\shade] ({2.05+1.55*\col},6.65) rectangle ++(1.55,1);
  \node[align=center] at ({2.825+1.55*\col},7.15) {$\atom$\\$f=\val$};
}
\foreach \col in {0,...,4} {
  \draw[steinRule] ({2.05+1.55*\col},6.65) -- ({2.05+1.55*\col},10.65);
}
\foreach \row in {0,...,4} {
  \draw[steinRule] (2.05,{6.65+\row}) -- (8.25,{6.65+\row});
}
\draw[steinInk!65,line width=.7pt] (2.05,7.65) -- (8.25,7.65);
\draw[line,draw=steinB] (8.43,7.15) -- (8.94,7.15);
\node[vertex,draw=steinB,fill=steinB!8] at (9.35,7.15) {$B$};
\draw[line,draw=steinB] (9.86,7.15) -- (10.61,7.15);
\node[rowval,text=steinB] at (11.35,7.15) {$2$};

\node[note] at (0,6.25)
  {Each cell is one point of $\{-1,1\}^4$, of probability $1/16$; its value is $J_2f=f\circ\ell_2$.};
\node[note] at (0,5.85)
  {The numbers of cells labelled $(A_1,A_2,B_1,B_2,B_3)$ are $(3,9,1,2,1)$.};
\end{tikzpicture}

\ifdefined\captionindent\setlength{\captionindent}{0pt}\fi
\small
\caption[Dyadic approximation and sampling in Step 2]{%
A numerical example of Step~2 with $n=1$ and $q=2$.
The atoms at levels $0$ and $1$ generate $\mathscr G_0$ and
$\mathscr G_1$, respectively. Averaging over the unused second block gives
\[
\begin{aligned}
 E_0^{[2]}f(A)
 &=\tfrac14\,0+\tfrac34\,4
   =\tfrac14(0+4+4+4)=3,\\
 E_0^{[2]}f(B)
 &=\tfrac14(-4)+\tfrac24\,2+\tfrac14\,8
   =\tfrac14(-4+2+2+8)=2.
\end{aligned}
\]
The first block distinguishes $--$, $-+$ and $+-$, but all three select
$A$ and give the same row. Conditional on any of these pairs, the unused
block is uniform, so the additional information does not change the
conditional average. Thus
\[
 \mathcal E_2^{(4)}J_2f=J_2E_0^{[2]}f,
 \qquad k_0=2,\quad k_1=4.
\]
Counting cells also gives the isometry: for every terminal $X$-valued
function $g$,
\[
 \|J_2g\|_{L^p(\{-1,1\}^4;X)}^p
 =\sum_{v\in\{A_1,A_2,B_1,B_2,B_3\}}
   \mu_2(v)\|g(v)\|^p
 =\|g\|_{L^p(\mu_2;X)}^p.
\]
These identities are exact for $\mu_2$. Under the original measure $\mu$,
the conditional averages of $f$ on $A$ and $B$ are $8/3$ and $13/5$.}
\label{fig:full-dyadic-Stein-sampling}
\end{minipage}
\end{figure}
\ifdefined\SteinFigureInManuscript
  \expandafter 
\fi

\endgroup

\emph{Step 2: Dyadic approximation and sampling.} After null atoms have been omitted, the atoms of \(\mathscr G_j\) form the vertices at level \(j\) of a finite tree. Each vertex at level \(j+1\) has as its parent the unique atom of \(\mathscr G_j\) containing it. We adjoin a root \(S\) at level \(-1\) and denote the set of terminal atoms by \(\mathscr T\). Functions measurable with respect to \(\mathscr G_n\) may then be identified, up to null sets, with functions on \(\mathscr T\).

For a child \(w\) of a vertex \(v\), the transition probability
\[
 \vartheta(w|v):=\frac{\mu(w)}{\mu(v)}
\]
is strictly positive, and \(\sum_{w\text{ child of }v}\vartheta(w|v)=1\). These transition probabilities can be approximated simultaneously by positive dyadic rational numbers with a common denominator. More precisely, for every sufficiently large integer \(q\) there are positive integers \(a_q(w|v)\) such that
\[
 \sum_{w\text{ child of }v}a_q(w|v)=2^q,
 \qquad
 \Bigl|2^{-q}a_q(w|v)-\vartheta(w|v)\Bigr|\leq 2^{-q}.
\]
To see this, suppose that \(v\) has \(d\) children. The integers \(\lfloor2^q\vartheta(w|v)\rfloor\) have sum \(2^q-r\) for some \(r\in\{0,\ldots,d-1\}\). Adding one to the rounded values for \(r\) distinct children gives the required sum and error bound. Since there are only finitely many edges and all transition probabilities are positive, all the rounded values are at least one once \(q\) is sufficiently large.

We write \(\vartheta_q(w|v):=2^{-q}a_q(w|v)\). For a vertex \(v\) at level \(j\), let \(S=v_{-1},v_0,\ldots,v_j=v\) be the branch joining it to the root, and define
\[
 \mu_q(S):=1,
 \qquad
 \mu_q(v):=\prod_{i=0}^j
                  \vartheta_q(v_i|v_{i-1}).
\]
The identity \(\sum_{w\text{ child of }v}\mu_q(w)=\mu_q(v)\) shows that the terminal weights define a probability measure on \(\mathscr T\). It also shows that, for every vertex \(v\),
\[
 \mu_q(v)=
 \sum_{\substack{a\in\mathscr T\\a\subseteq v}}\mu_q(a).
\]
All these weights are positive. Moreover, \(\mu_q(v)\to\mu(v)\) at every vertex, since the corresponding product of the original transition probabilities telescopes to \(\mu(v)\).

For \(a\in\mathscr T\), let \(v_j(a)\) be its ancestor at level \(j\). For \(f:\mathscr T\to X\), conditional expectation for \(\mu_q\) onto the partition at level \(j\) is therefore given by
\[
 (E_j^{[q]}f)(a)
 =
 \sum_{\substack{b\in\mathscr T\\b\subseteq v_j(a)}}
       \frac{\mu_q(b)}{\mu_q(v_j(a))}\,f(b).
\]

We next realise these conditional expectations on a dyadic cube. For each vertex \(v\), there is a partition
\[
 \{1,\ldots,2^q\}
 =
 \mathop{\dot\bigcup}_{w\text{ child of }v} I_q(w|v),
 \qquad
 |I_q(w|v)|=a_q(w|v).
\]
Set \(d:=(n+1)q\). A fixed bijection between \(\D^q\) and \(\{1,\ldots,2^q\}\) turns the successive blocks of \(q\) coordinates of \(\D^d\) into independent uniform random variables \(U_0,\ldots,U_n\) with values in \(\{1,\ldots,2^q\}\).

These variables determine vertices \(V_{-1},V_0,\ldots,V_n\) recursively. The initial vertex is \(V_{-1}=S\), and, when \(V_{j-1}=v\), the next vertex is the unique child \(w\) for which \(U_j\in I_q(w|v)\). Independence of the blocks gives
\[
 \mathbb P(V_j=w\,|\,V_{j-1}=v)
 =\frac{a_q(w|v)}{2^q}
 =\vartheta_q(w|v).
\]
It follows by induction that \(\mathbb P(V_j=v)=\mu_q(v)\). In particular, the terminal map \(\ell_q:=V_n\) has law \(\mu_q\).

The linear map
\[
 J_qf:=f\circ\ell_q
\]
is consequently an isometry from \(L^p(\mathscr T,\mu_q;X)\) into \(L^p(\D^d;X)\), since
\[
 \|J_qf\|_{L^p(\D^d;X)}^p
 =
 \sum_{a\in\mathscr T}\mu_q(a)\|f(a)\|^p
 =
 \|f\|_{L^p(\mathscr T,\mu_q;X)}^p.
\]
To verify the relation between conditional expectations, let \(k_j:=(j+1)q\) and write points of \(\D^d\) as \((\alpha,z)\) with \(\alpha\in\D^{k_j}\) and \(z\in\D^{d-k_j}\). The coordinates in \(\alpha\) form the first \(j+1\) blocks and determine a vertex \(v\) at level \(j\). For each terminal atom \(a\subseteq v\), the number of choices of \(z\) for which \(\ell_q(\alpha,z)=a\) is the product of the integers \(a_q(w|u)\) along the edges from \(v\) to \(a\). Hence
\[
 2^{-(d-k_j)}
 \#\bigl\{z\in\D^{d-k_j}:\ \ell_q(\alpha,z)=a\bigr\}
 =
 \begin{cases}
 \displaystyle\frac{\mu_q(a)}{\mu_q(v)},
      & a\subseteq v,\\[6pt]
 0,   & a\not\subseteq v.
 \end{cases}
\]
At \(j=n\) the same formula holds with \(z\) the empty tuple. Averaging \((J_qf)(\alpha,z)\) over \(z\in\D^{d-k_j}\) therefore gives
\[
 \bigl(\mathcal E_{k_j}^{(d)}J_qf\bigr)(\alpha,z)
 =
 \sum_{\substack{a\in\mathscr T\\a\subseteq v}}
       \frac{\mu_q(a)}{\mu_q(v)}\,f(a)
 =
 \bigl(J_qE_j^{[q]}f\bigr)(\alpha,z).
\]
Thus
\[
 \mathcal E_{k_j}^{(d)}J_q=J_qE_j^{[q]},
 \qquad 0\leq j\leq n.
\]
In particular, the conditional average depends on the first \(k_j\) coordinates only through the vertex \(v\) they determine. Figure~\ref{fig:full-dyadic-Stein-sampling} illustrates this construction.

\smallskip
\emph{Step 3: Passing to the original probabilities.} The dyadic Stein estimate, applied to the functions \(J_qf_\ell\) at the indices \(k_{\tau(\ell)}\), gives
\[
 \Bigl\|\sum_{\ell=1}^m\varepsilon_\ell
              E_{\tau(\ell)}^{[q]}f_\ell
       \Bigr\|_{\Rad_m(L^p(\mu_q;X))}
 \leq C
 \Bigl\|\sum_{\ell=1}^m\varepsilon_\ell f_\ell
       \Bigr\|_{\Rad_m(L^p(\mu_q;X))}.
\]
Here we have used the intertwining identity and the fact that the linear isometry \(J_q\) preserves the \(L^p\)-norm of each realisation of a random sum.

For every terminal atom \(a\) and every \(j\), the explicit formula for \(E_j^{[q]}\) gives
\[
 (E_j^{[q]}f_\ell)(a)
 \to
 \sum_{\substack{b\in\mathscr T\\b\subseteq v_j(a)}}
       \frac{\mu(b)}{\mu(v_j(a))}\,f_\ell(b)
 =
 \E(f_\ell|\mathscr G_j)(a).
\]
Indeed, every numerator converges and the limiting denominator \(\mu(v_j(a))\) is strictly positive. Since \(\mathscr T\) is finite and \(\mu_q(a)\to\mu(a)\), it follows that, for each fixed realisation of the external variables,
\[
 \Bigl\|\sum_{\ell=1}^m\varepsilon_\ell
                 E_{\tau(\ell)}^{[q]}f_\ell
       \Bigr\|_{L^p(\mu_q;X)}
 \to
 \Bigl\|\sum_{\ell=1}^m\varepsilon_\ell
                 \E(f_\ell|\mathscr G_{\tau(\ell)})
       \Bigr\|_{L^p(\mu;X)}.
\]
The corresponding convergence
\[
 \Bigl\|\sum_{\ell=1}^m\varepsilon_\ell f_\ell\Bigr\|_{L^p(\mu_q;X)}
 \to
 \Bigl\|\sum_{\ell=1}^m\varepsilon_\ell f_\ell\Bigr\|_{L^p(\mu;X)}
 \qquad (q\to\infty)
\]
follows in the same way.

All these norms are bounded by
\[
 M:=\sum_{\ell=1}^m\max_{a\in\mathscr T}\|f_\ell(a)\|,
\]
independently of \(q\) and of the external variables. For the sums involving conditional expectations, the bound follows because each value of \(E_j^{[q]}f_\ell\) and \(\E(f_\ell|\mathscr G_j)\) is a convex combination of the values \(f_\ell(a)\), \(a\in\mathscr T\). Dominated convergence, applied to the squares of the displayed norms, therefore gives convergence of both randomisation norms. Passing to the limit proves the finite-tree Stein estimate with constant \(C\).

By Step~1, we obtain \(\mathfrak s_{p,X}\leq C\). Together with the reverse inequality \(\mathfrak s^{\rm dyad}_{p,X}\leq\mathfrak s_{p,X}\), this proves equality of the optimal constants.

The same argument applies when the second moment is replaced by any fixed \(r\)th moment, \(1\leq r<\infty\), and when signs are used instead of Steinhaus variables. The contractions and isometries act on each realisation of the random sums, and the bound by \(M\) permits dominated convergence for the \(r\)th powers as well. For each of these conventions, the full and dyadic optimal constants are therefore equal.
\end{proof}

\subsection{Independence of the exponent}

We show next that if a Banach space has the \(p\)-dyadic Stein property for some \(1<p<\infty\), then it has the \(p\)-dyadic Stein property for all \(1<p<\infty\).
By Theorems~\ref{thm:full-dyadic-Stein}, the same is then true for the \(p\)-Stein property. This will allow us to speak of the \emph{(dyadic) Stein property} without reference to a particular chose of exponent \(1<p<\infty\).

\begin{theorem}[Independence of the exponent]
\label{thm:dyadic-Stein-p-independent}
Let \(X\) be a Banach space and let \(1<p,q<\infty\). If \(X\) has the \(p\)-dyadic Stein property, then it has the \(q\)-dyadic Stein property. Moreover, \(\mathfrak s^{\rm dyad}_{q,X}\) is bounded in terms of \(p\), \(q\), and \(\mathfrak s^{\rm dyad}_{p,X}\) only.
\end{theorem}

\begin{proof}
We first compare the Stein constant with the norms of suitable diagonal operators on \(L^r(\D^d;\Rad_m(X))\). We then extrapolate these operator bounds to the required exponent. The extrapolation uses a dyadic Calder\'on--Zygmund decomposition and duality, as in the classical extrapolation of martingale transforms \cite[Chapter 4]{HNVW1}, \cite[Chapter~6]{Pisier-Martingales}.

\smallskip
\emph{Step 1: Comparing \(\mathfrak s^{\rm dyad}_{r,X}\) with diagonal operator norms.}

For \(m\geq1\) put \(Y_m:=\Rad_m(X)\), and let \(Q_1,\ldots,Q_m\) denote the canonical coordinate projections on \(Y_m\). Fix \(d,m\geq1\) and a map \(\tau:\{1,\ldots,m\}\to\{0,\ldots,d\}\).

For \(y=\sum_{\ell=1}^m\varepsilon_\ell x_\ell\in Y_m\), the coefficient maps are
\[
 C_\ell y:=\E_\varepsilon(\overline{\varepsilon_\ell}y)=x_\ell.
\]
The Cauchy--Schwarz inequality gives \(\|C_\ell y\|_X\leq\|y\|_{Y_m}\). Thus every \(F\in L^r(\D^d;Y_m)\) has the unique representation
\[
 F=\sum_{\ell=1}^m\varepsilon_\ell f_\ell,
 \qquad f_\ell=C_\ell F\in L^r(\D^d;X).
\]
In particular, \(Q_\ell y=\varepsilon_\ell C_\ell y\) is a bounded projection. We define
\[
 P_\tau F:=\sum_{\ell=1}^m
       \varepsilon_\ell\mathcal E_{\tau(\ell)}^{(d)}f_\ell
       =\sum_{\ell=1}^m Q_\ell\mathcal E_{\tau(\ell)}^{(d)}F.
\]
Here conditional expectation acts on the cube variable and \(Q_\ell\) acts on the value in \(Y_m\). They commute because bounded linear maps commute with the finite averages defining these conditional expectations.

For clarity, the comparison of the two mixed norms is
\begin{align*}
 \|F\|_{L^r(\D^d;Y_m)}
 &=\Bigl(\E_x\bigl(\E_\varepsilon
       \|F(x,\varepsilon)\|_X^2\bigr)^{r/2}\Bigr)^{1/r}\\
 &\eqsim_r\Bigl(\E_x\E_\varepsilon
       \|F(x,\varepsilon)\|_X^r\Bigr)^{1/r}\\
 &=\Bigl(\E_\varepsilon
       \|F(\cdot,\varepsilon)\|_{L^r(\D^d;X)}^r\Bigr)^{1/r}\\
 &\eqsim_r\Bigl(\E_\varepsilon
       \|F(\cdot,\varepsilon)\|_{L^r(\D^d;X)}^2\Bigr)^{1/2}.
\end{align*}
The first comparison is the Kahane--Khintchine inequality in \(X\), applied for each \(x\), and the second is the same inequality in \(L^r(\D^d;X)\). Fubini's theorem gives the intervening equality. The constants depend only on \(r\), and the same comparisons apply to \(P_\tau F\). As \(\tau\) varies, the list \(\mathcal E_{\tau(1)}^{(d)},\ldots,\mathcal E_{\tau(m)}^{(d)}\) ranges over all lists of \(m\) operators from \(\{\mathcal E_0^{(d)},\ldots,\mathcal E_d^{(d)}\}\), with repetitions allowed. The functions \(f_1,\ldots,f_m\) may be chosen arbitrarily in \(L^r(\D^d;X)\).

Taking the supremum over \(d,m\geq1\) and all maps \(\tau\) as above, we obtain
\begin{align}\label{eq:Stein-adapted-projection} \sup_{d,m,\tau}
\|P_\tau\|_{\calL(L^r(\D^d;Y_m))} \,\eqsim_r\, \mathfrak s^{\rm dyad}_{r,X}.
\end{align}

\smallskip
\emph{Step 2: Downward extrapolation.} Suppose now that for some \(1<p< \infty\) we have \(\mathfrak s^{\rm dyad}_{p,X}<\infty\). In this step we prove that if \(1<q<p\), then \(\mathfrak s^{\rm dyad}_{q,X}<\infty\). We do this by applying \eqref{eq:Stein-adapted-projection} with exponents \(p\) and \(q\).

Throughout the remaining steps, we write
\[
 C:=\sup_{d,m,\tau}
 \|P_\tau\|_{\calL(L^p(\D^d;Y_m))}.
\]
By \eqref{eq:Stein-adapted-projection}, this constant is finite and is bounded in terms of \(p\) and \(\mathfrak s^{\rm dyad}_{p,X}\) only.

Fix \(d,m\geq1\) and a map \(\tau\) as above, and let \(F\in L^1(\D^d;Y_m)\). Let \(\lambda>0\). If \(\lambda<\|F\|_{1}\), then
\begin{align}\label{eq:weak11-1}
 \big|\{\|P_\tau F\|_{Y_m}>\lambda\}\big|
 \leq1
 <\frac{\|F\|_{1}}{\lambda}.
\end{align}

Suppose next that \(\lambda\geq\|F\|_{1}\). To prove a similar weak type \((1,1)\) estimate, we use a Calder\'on--Zygmund decomposition argument.

We consider the atoms of all the \(\sigma\)-algebras in the canonical filtration on \(\D^d\). These are the sets
\[
 A(a_1,\ldots,a_k)
 :=\{x\in\D^d:\,x_1=a_1,\ldots,x_k=a_k\},
\]
with \(0\leq k\leq d\) and \((a_1,\ldots,a_k)\in\D^k\). An atom at level \(k\) has probability \(2^{-k}\); if \(k\geq1\), its parent is the atom at level \(k-1\) obtained by omitting the condition on the \(k\)-th coordinate.

For each atom \(A\), we write
\[
 \langle F\rangle_A:=\frac1{|A|}\int_A F,
 \qquad
 \langle\|F\|_{Y_m}\rangle_A
 :=\frac1{|A|}\int_A\|F\|_{Y_m},
\]
where \(|A| = \#A/2^d\) denotes the uniform measure of \(A\). Let \(\mathcal A\) be the collection of atoms which are maximal with respect to inclusion among those satisfying
\[
 \langle\|F\|_{Y_m}\rangle_A>\lambda.
\]
Two atoms of the canonical filtration are either disjoint or one contains the other, so the members of \(\mathcal A\) are pairwise disjoint. Moreover, the assumption \(\lambda\geq\|F\|_1\) ensures that \(\D^d\) does not belong to \(\mathcal A\). Each \(A\in\mathcal A\) therefore has an immediate parent \(A^+\), and maximality gives
\[
 \langle\|F\|_{Y_m}\rangle_{A^+}\leq\lambda.
\]
Since \(|A^+|=2|A|\), it follows that
\[
 \|\langle F\rangle_A\|_{Y_m}
 \leq \langle\|F\|_{Y_m}\rangle_A
 \leq 2\langle\|F\|_{Y_m}\rangle_{A^+}
 \leq 2\lambda.
\]

We now put
\[
 B:=\bigcup_{A\in\mathcal A}A,
 \qquad
 b_A:={\bf1}_A(F-\langle F\rangle_A),
 \qquad
 g:=F-\sum_{A\in\mathcal A}b_A.
\]
Thus \(g=F\) on \(\D^d\setminus B\), whereas \(g=\langle F\rangle_A\) on each \(A\in\mathcal A\). Outside \(B\), we have \(\|F\|_{Y_m}\leq\lambda\): otherwise the singleton containing such a point would satisfy the defining inequality and would be contained in a member of \(\mathcal A\). Consequently,
\begin{equation}
\label{eq:Stein-CZ-good}
 \|g\|_{L^\infty}\leq2\lambda,
 \qquad
 \|g\|_1\leq\|F\|_1,
 \qquad
 |B|\leq\frac{\|F\|_1}{\lambda}.
\end{equation}
The second estimate follows by summing \(|A|\|\langle F\rangle_A\|_{Y_m} \leq\int_A\|F\|_{Y_m}\) over \(A\in\mathcal A\) and adding the integral over \(\D^d\setminus B\). The third follows by summing \(\lambda|A|<\int_A\|F\|_{Y_m}\) over \(A\in\mathcal A\).

Each function \(b_A\) is supported in \(A\) and has integral zero:
\[
 \int_A b_A=0.
\]
These two properties imply that \(P_\tau b_A\) is also supported in \(A\). Indeed, suppose that \(A\) is an atom at level \(k\). If \(\tau(\ell)<k\), then every atom at level \(\tau(\ell)\) either contains \(A\) or is disjoint from \(A\). The support and zero-integral properties therefore give
\[
 \mathcal E_{\tau(\ell)}^{(d)}b_A=0.
\]
If \(\tau(\ell)\geq k\), then \(A\) is measurable with respect to the first \(\tau(\ell)\) coordinates, so
\[
 \mathcal E_{\tau(\ell)}^{(d)}b_A
 ={\bf1}_A\mathcal E_{\tau(\ell)}^{(d)}b_A.
\]
Since each \(Q_\ell\) acts pointwise on the values in \(Y_m\), it also preserves support. Hence
\begin{equation*}
 \operatorname{supp}P_\tau b_A\subseteq A.
\end{equation*}
Applying \(P_\tau\) to the decomposition \(F=g+\sum_{A\in\mathcal A}b_A\), we conclude that \(P_\tau F=P_\tau g\) on \(\D^d\setminus B\).

By Chebyshev's inequality, the \(L^p\)-bound, and \eqref{eq:Stein-CZ-good},
\begin{equation}\label{eq:weak11-2}
\begin{aligned}
 \big|\{\|P_\tau F\|_{Y_m}>\lambda\}\big|
 &\leq \frac{\|F\|_1}{\lambda}
       +\lambda^{-p}\|P_\tau g\|_{L^p}^p\\
 &\leq \frac{\|F\|_1}{\lambda}
       +C^p\lambda^{-p}\|g\|_{L^\infty}^{p-1}\|g\|_1
 \leq (1+2^{p-1}C^p)\frac{\|F\|_1}{\lambda},
\end{aligned}
\end{equation}

Combining \eqref{eq:weak11-1} and \eqref{eq:weak11-2}, we obtain the uniform weak type \((1,1)\) constant \(A_p:=1+2^{p-1}C^p\). We give the interpolation step explicitly. Fix \(1<q<p\) and, for each \(\lambda>0\), decompose
\[
 F=F_{>\lambda}+F_{\leq\lambda},
 \qquad
 F_{>\lambda}:={\bf1}_{\{\|F\|_{Y_m}>\lambda\}}F.
\]
By linearity and the triangle inequality, at every point where \(\|P_\tau F\|_{Y_m}>\lambda\), at least one of \(\|P_\tau F_{>\lambda}\|_{Y_m}\) and \(\|P_\tau F_{\leq\lambda}\|_{Y_m}\) exceeds \(\lambda/2\). Applying the weak type estimate to \(F_{>\lambda}\) and the \(L^p\) estimate to \(F_{\leq\lambda}\) therefore gives
\begin{align*}
 |\{\|P_\tau F\|>\lambda\}|
 \leq\frac{2A_p}{\lambda}
       \int_{\{\|F\|>\lambda\}}\|F\|
 +\frac{2^pC^p}{\lambda^p}
       \int_{\{\|F\|\leq\lambda\}}\|F\|^p.
\end{align*}
We multiply by \(q\lambda^{q-1}\) and integrate over \(\lambda>0\). The distribution formula for the \(L^q\)-norm and Fubini's theorem then give
\begin{align*}
 \|P_\tau F\|_q^q
 &\leq 2qA_p\!\int_{\D^d} \|F(x)\|
          \int_0^{\|F(x)\|}\!\!\lambda^{q-2}\dd\lambda\dd x\\
 &\quad+2^pqC^p\!\!\int_{\D^d}\|F(x)\|^p
          \int_{\|F(x)\|}^{\infty}\!\lambda^{q-p-1}\dd\lambda\dd x
 =\Bigl(\frac{2qA_p}{q-1}
          +\frac{2^pqC^p}{p-q}\Bigr)\|F\|_q^q.
\end{align*}
At points where \(F=0\), the corresponding integrands are zero. The restrictions \(q>1\) and \(q<p\) ensure convergence of the two inner integrals. Thus
\begin{equation}
\label{eq:Stein-extrapolate-down}
 \sup_{d,m,\tau}\|P_\tau\|_{\calL(L^q(\D^d;Y_m))}
 \leq\Bigl(\frac{2q(1+2^{p-1}C^p)}{q-1}
          +\frac{2^pqC^p}{p-q}\Bigr)^{1/q}<\infty.
\end{equation}
Together with \eqref{eq:Stein-adapted-projection}, this proves downward extrapolation with the asserted dependence of the constant.

\smallskip
\emph{Step 3: Upward extrapolation by duality.} Suppose next that for some \(1<p< \infty\) we have \(\mathfrak s^{\rm dyad}_{p,X}<\infty\). In this step we prove that if \(p<q<\infty\), then \(\mathfrak s^{\rm dyad}_{q,X}<\infty\). For this we combine the result of Step 2 with a duality argument.

We will use the actual dual space \(Y_m^*\), without identifying it with a space of random sums over \(X^*\). Since \(\D^d\) is a finite atomic measure space, \cite[Proposition~1.3.3]{HNVW1} gives the isometric identification
\[
 L^p(\D^d;Y_m)^*=L^{p'}(\D^d;Y_m^*).
\]
For functions \(F\) and \(G\) on \(\D^d\) with values in \(Y_m\) and \(Y_m^*\), respectively, averaging on the atoms gives
\[
 \int_{\D^d}\langle Q_\ell\mathcal E_{\tau(\ell)}^{(d)}F,G\rangle
 =\int_{\D^d}\langle F,
       \mathcal E_{\tau(\ell)}^{(d)}Q_\ell^*G\rangle.
\]
The operator \(Q_\ell^*\) commutes with conditional expectation for the same reason as \(Q_\ell\). Consequently, the adjoint of \(P_\tau\) is given on \(L^{p'}(\D^d;Y_m^*)\) by
\[
 P_\tau^* = \sum_{\ell=1}^mQ_\ell^*\mathcal E_{\tau(\ell)}^{(d)}.
\]
Since an operator and its adjoint have the same norm, we have
\[
 \|P_\tau^*\|_{\calL(L^{p'}(\D^d;Y_m^*))}\leq C.
\]
To see that the same argument applies to this adjoint, suppose that \(b\) is \(Y_m^*\)-valued, supported in a level-\(k\) atom \(A\), and satisfies \(\int_A b=0\). For \(\tau(\ell)<k\), the conditional expectation of \(b\) vanishes, whereas for \(\tau(\ell)\geq k\) it is supported in \(A\). The pointwise operator \(Q_\ell^*\) preserves that support. Hence \(P_\tau^*b\) is supported in \(A\). The decomposition and interpolation in Step~2 used only this support property and a strong bound at the starting exponent; they did not use any further structure of \(Y_m\). They therefore apply with \(Y_m^*\) in place of \(Y_m\) and \(p'\) in place of \(p\). If \(q>p\), then \(1<q'<p'\), so they give
\[
 \sup_{d,m,\tau}
 \|P_\tau^*\|_{\calL(L^{q'}(\D^d;Y_m^*))}
 \leq C_{p,q,C}<\infty.
\]
This constant depends only on \(p,q,C\). The displayed formula for \(P_\tau^*\) is independent of the exponent. Using the finite-atomic duality once more, we obtain
\[
 \|P_\tau F\|_{L^q(Y_m)}
 =\sup_{\|G\|_{L^{q'}(Y_m^*)}\leq1}
   |\langle F,P_\tau^*G\rangle|
 \leq\|P_\tau^*\|_{L^{q'}(Y_m^*)\to L^{q'}(Y_m^*)}
       \|F\|_{L^q(Y_m)}.
\]
Consequently,
\[
 \sup_{d,m,\tau} \|P_\tau\|_{\calL(L^q(\D^d;Y_m))}<\infty, \qquad p<q<\infty.
\]
Together with \eqref{eq:Stein-extrapolate-down} and \eqref{eq:Stein-adapted-projection}, this proves the result.
\end{proof}

\begin{remark}[An explicit extrapolation bound]
Let \(K_{a,b}\) be the Kahane--Khintchine constant for the complex scalar field as defined in \cite[Theorem~6.2.4]{HNVW2}, so that \[\|Z\|_{L^a}\leq K_{a,b}\|Z\|_{L^b}\] for finite Steinhaus sums \(Z\). Put \(\kappa_r:=K_{\max\{r,2\},\min\{r,2\}}\) and, for \(1<b<a<\infty\), define
\[
 D(a,b):=\left(\frac{2b(1+2^{a-1})}{b-1}+\frac{2^a b}{a-b}\right)^{1/b}.
\]
Then we have the explicit bound
\[
 \mathfrak s_{q,X}\leq
 \begin{cases}
  \kappa_p\kappa_q D(p,q)\,\mathfrak s_{p,X}, & q<p,\\
  \mathfrak s_{p,X}, & q=p,\\
  \kappa_p\kappa_q D(p',q')\,\mathfrak s_{p,X}, & q>p,
 \end{cases}
\]
where \(p'\) and \(q'\) are the conjugate exponents. The same bounds hold for the dyadic constants.

Indeed, in Step~1 write \[A_r(F):=\|F\|_{\Rad_m(L^r(\D^d;X))}, \qquad B_r(F):=\|F\|_{L^r(\D^d;Y_m)}.\] Minkowski's integral inequality gives \(B_r(F)\leq A_r(F)\) for \(r\geq2\) and \(A_r(F)\leq B_r(F)\) for \(r\leq2\). Monotonicity of the probability norms and the Kahane--Khintchine inequality at each cube point place the joint \(L^r\)-norm between \(A_r(F)\) and \(\kappa_r B_r(F)\) for \(r\geq2\), and between \(\kappa_r^{-1} B_r(F)\) and \(A_r(F)\) for \(r\leq2\). Thus the operator norms with respect to \(A_r\) and \(B_r\) are bounded by \(\kappa_r\) times one another.

In particular, the uniform \(L^p\)-bound \(C\) of Step~2 satisfies \(C\leq\kappa_p\mathfrak s_{p,X}^{\rm dyad}\). For \(C>0\), the operators \(P_\tau/C\) and \(P_\tau^*/C\) have the support property used in Step~2 and norm at most one at exponents \(p\) and \(p'\), respectively. Applying \eqref{eq:Stein-extrapolate-down} with unit starting bound and rescaling gives the uniform \(L^q\)-bound \(C D(p,q)\) when \(q<p\), and \(C D(p',q')\) when \(q>p\) by duality as in Step~3. The case \(C=0\) is immediate. The norm comparison at \(q\) gives the dyadic estimate, and Theorem~\ref{thm:full-dyadic-Stein} identifies the dyadic and full constants.
\end{remark}

\begin{remark}
The proof gives the uniform weak type \((1,1)\) estimate
\[
 \lambda\big|\{\|P_\tau F\|_{Y_m}>\lambda\}\big|
 \leq(1+2^{p-1}C^p)\|F\|_1
\]
for the operators \(P_\tau\). The corresponding uniform strong \(p=1\) estimate fails already for the scalar field, as in the classical endpoint failure of Stein's inequality \cite[p.~105]{Stein-LP}. A dyadic example is obtained by setting
\[
 A_j=\{x\in\D^m:x_1=\cdots=x_j=1\},
 \qquad 0\leq j\leq m,
\]
with \(A_0=\D^m\), and taking \(f_j={\bf1}_{A_j\setminus A_{j+1}}\) for \(0\leq j<m\). These supports are disjoint, so
\[
 \|(f_j)_{j=0}^{m-1}\|_{L^\infty(\D^m;\ell^2_m)}=1.
\]
At the point \((1,\ldots,1)\), the conditional average \(\mathcal E_j^{(m)}f_j\) equals \(1/2\): among the continuations of the first \(j\) signs, exactly half have \(x_{j+1}=-1\). Thus the diagonal operator
\[
 (f_j)_{j=0}^{m-1}\mapsto
 (\mathcal E_j^{(m)}f_j)_{j=0}^{m-1}
\]
has norm at least \(\sqrt m/2\) on \(L^\infty(\D^m;\ell^2_m)\). Its adjoint under the integral pairing is the same diagonal operator, so its norm on \(L^1(\D^m;\ell^2_m)\) has the same lower bound. Finally, \(\Rad_m(\mathbb C)=\ell^2_m\) isometrically, and the mixed-norm comparison in Step~1 remains valid for \(r=1\), with constants independent of \(m\). The scalar dyadic conditional expectations therefore cannot have a uniform \(R\)-bound on \(L^1\).
\end{remark}

\subsection{Self-duality}

Self-duality of the Stein property is implicit in Eskenazis's argument \cite[proof of Proposition~5, (37)--(41)]{Eskenazis}. By \cite[Lemma~4]{Eskenazis}, the Stein property implies \(K\)-convexity, which gives a uniform identification of the dual of a space of random sums with the corresponding space over \(X^*\). The general duality theorem for \(R\)-bounded families on spaces with non-trivial type is given in \cite[Theorem~5.2]{KwapienVeraarWeis}.
For the reader's convenience we include a direct argument.

\begin{proposition}[Self-duality of the Stein property]
\label{prop:Stein-selfdual}
A Banach space \(X\) has the Stein property if and only if \(X^*\) does.
\end{proposition}

\begin{proof}
Suppose that \(X\) has the Stein property. Theorems \ref{thm:full-dyadic-Stein} and \ref{thm:dyadic-Stein-p-independent} allow us to work at exponent \(p=2\). As we have already noted, the Stein property implies that \(X\) is \(K\)-convex.

In this proof we use randomisation by Rademacher variables \(r_j\). The comparison between Rademacher and Steinhaus sums has constants independent of the number of summands and of the Banach space \cite[Section~6.1.a]{HNVW2}; hence either randomisation gives the same qualitative Stein property. The sign convention also allows all duality arguments below to take place on finite cubes.

We write \(\Rad_m^r(X)\) for the space of Rademacher sums \(\sum_{j=1}^m r_jx_j\) with its \(L^2(\D^m;X)\)-norm. The assocated {\em Rademacher projection} is defined by
\[
 (\Pi_m^XF)(r):=\sum_{j=1}^m r_j\E_s(s_jF(s)).
\]
Here \(s\) is an independent uniform sign vector. The identities \(\E_s(s_js_k)=\delta_{jk}\) show that \(\Pi_m^X\) is a projection onto \(\Rad_m^r(X)\). By \(K\)-convexity on \(X\) \cite{Pisier}, see also \cite[Proposition 7.4.5]{HNVW2},
\[
 K_X:=\sup_{m\geq1}
 \|\Pi_m^X\|_{\calL(L^2(\D^m;X))}<\infty.
\]

We first establish the required duality for random sums. For \(g=\sum_jr_jx_j^*\in\Rad_m^r(X^*)\), define
\[
 (J_mg)\Bigl(\sum_jr_jx_j\Bigr)
 :=\sum_j\langle x_j,x_j^*\rangle.
\]
Orthogonality of the signs and the Cauchy--Schwarz inequality give
\[
 |(J_mg)(h)|
 =|\E_r\langle h(r),g(r)\rangle|
 \leq\|h\|_{L^2(\D^m;X)}\|g\|_{L^2(\D^m;X^*)}.
\]
Thus \(\|J_m\|\leq1\). If \(J_mg=0\), then \((J_mg)(r_jx)=\langle x,x_j^*\rangle=0\) for every \(j\) and every \(x\in X\). Hence \(x_j^*=0\) for all \(j\), which proves that \(J_m\) is injective.

For surjectivity, let \(\psi\in\Rad_m^r(X)^*\). The Hahn--Banach theorem extends \(\psi\) to a functional on \(L^2(\D^m;X)\) with the same norm. Since the cube is finite, its dual is \(L^2(\D^m;X^*)\) under the integral pairing \cite[Proposition~1.3.3]{HNVW1}. There is therefore a \(G\in L^2(\D^m;X^*)\) representing this extension with \(\|G\|_2=\|\psi\|\). Directly from the formula for the projection,
\[
 \E_r\langle\Pi_m^XF(r),G(r)\rangle
 =\E_r\langle F(r),\Pi_m^{X^*}G(r)\rangle.
\]
Thus \(\Pi_m^{X^*}\) is the adjoint of \(\Pi_m^X\) and has norm at most \(K_X\). Set \(g=\Pi_m^{X^*}G\). For \(h\in\Rad_m^r(X)\),
\[
 (J_mg)(h)
 =\E_r\langle\Pi_m^Xh(r),G(r)\rangle
 =\psi(h),
 \qquad \|g\|_2\leq K_X\|\psi\|.
\]
We have proved that \(J_m\) is an isomorphism and that \(\|J_m^{-1}\|\leq K_X\), uniformly in \(m\).

Now fix \(d,m\geq1\) and \(k_j\in\{0,\ldots,d\}\). On \(Z_X:=L^2(\D^d;\Rad_m^r(X))\) define
\[
 P_X\Bigl(\sum_jr_jf_j\Bigr)
 :=\sum_jr_j\mathcal E_{k_j}^{(d)}f_j,
\]
and define \(Z_{X^*}\) and \(P_{X^*}\) in the same way. Fubini's theorem identifies the \(Z_X\)-norm exactly with the external \(L^2\) randomisation norm on \(L^2(\D^d;X)\). Consequently, the dyadic Stein constant for Rademacher randomisation at \(p=2\) bounds \(\|P_X\|\) independently of \(d\), \(m\) and the chosen levels \(k_1,\ldots,k_m\).

Pointwise application of \(J_m\), followed by finite-cube duality, gives an isomorphism
\[
 J:Z_{X^*}\to Z_X^*,
 \qquad \|J\|\leq1,\quad\|J^{-1}\|\leq K_X.
\]
For \(F=\sum_jr_jf_j\) and \(G=\sum_jr_jg_j\), self-adjointness of the conditional averages gives
\begin{align*}
 (JP_{X^*}G)(F)
 =\sum_j\int_{\D^d}
       \langle f_j,\mathcal E_{k_j}^{(d)}g_j\rangle
 =\sum_j\int_{\D^d}
       \langle\mathcal E_{k_j}^{(d)}f_j,g_j\rangle
 =(JG)(P_XF).
\end{align*}
Hence \(JP_{X^*}=P_X^*J\), and therefore
\[
 \|P_{X^*}\|
 \leq\|J^{-1}\|\,\|P_X^*\|\,\|J\|
 \leq K_X\|P_X\|.
\]
The bound for \(P_{X^*}\) holds for every \(d,m\geq1\) and every list \(k_1,\ldots,k_m\in\{0,\ldots,d\}\), with repetitions allowed. By the same application of Fubini's theorem, this is the dyadic Stein inequality for \(X^*\) with Rademacher randomisation. The comparison of randomisations and Theorem~\ref{thm:full-dyadic-Stein} give its full Stein property.

Conversely, if \(X^*\) has the Stein property, the implication just proved gives it for \(X^{**}\). The canonical map \(i:X\to X^{**}\) is an isometry, and for every Bochner integrable \(X\)-valued \(f\) we have
\[
 \E(i\circ f\mid\mathscr G)=i\circ\E(f\mid\mathscr G).
\]
This identity holds first for simple functions and then by continuity of conditional expectation. The Stein inequality in \(X^{**}\), applied to the functions \(i\circ f_j\), is thus exactly the same inequality in \(X\). This completes the proof.
\end{proof}

It remains an open question whether the Stein property implies UMD for arbitrary Banach spaces. Asking only for \(\mathrm{UMD}^+\) appears weaker, but self-duality makes the two general assertions equivalent.

\begin{corollary}
\label{cor:Stein-UMD}
The following assertions are equivalent:
\begin{enumerate}[\rm(1)]
\item\label{it:Stein-UMD1} Every Banach space with the Stein property is UMD\(^+\).
\item\label{it:Stein-UMD2} Every Banach space with the Stein property is UMD.
\end{enumerate}
If these assertions hold, then every UMD\(^+\) space is UMD.
\end{corollary}

\begin{proof}
It is clear from the equivalent randomised definition of the UMD property \cite[Proposition 4.2.3]{HNVW1} that \ref{it:Stein-UMD2} implies \ref{it:Stein-UMD1}. For the converse, suppose that \ref{it:Stein-UMD1} holds and that \(X\) has the Stein property. Proposition~\ref{prop:Stein-selfdual} gives the Stein property for \(X^*\), so the hypothesis makes both \(X\) and \(X^*\) UMD\(^+\).

By a result of Garling \cite{Garling}, the UMD\(^+\) property is independent of the exponent \(1<p<\infty\) occurring in the definition; see also \cite[Section~2.4]{KaltonLoristWeis}. As a result, we may use exponent \(p=2\). There are constants \(C_X,C_{X^*}\) such that, for \(Y=X\) or \(X^*\), every \(d\geq1\), and every finite \(Y\)-valued martingale difference sequence \((h_j)_{j=1}^m\) on \(\D^d\),
\[
 \Bigl(\E_r\Bigl\|\sum_{j=1}^mr_jh_j\Bigr\|_{L^2(\D^d;Y)}^2\Bigr)^{1/2}
 \leq C_Y\Bigl\|\sum_{j=1}^mh_j\Bigr\|_{L^2(\D^d;Y)}.
\]
We show how the upper estimate in \(X^*\) supplies the reverse estimate in \(X\).

On \(\D^d\), write \(\mathcal D_j=\mathcal E_j-\mathcal E_{j-1}\) for \(1\leq j\leq d\). If \((h_j)_{j=1}^d\) is a martingale difference sequence for the canonical filtration, then \(\mathcal D_jh_j=h_j\). For \(g\in L^2(\D^d;X^*)\), orthogonality of the external signs and self-adjointness of the conditional expectations give
\begin{align*}
 \E_r\int_{\D^d}\Bigl\langle\sum_{j=1}^dr_jh_j,
                          \sum_{k=1}^d r_k\mathcal D_kg\Bigr\rangle
 &=\sum_{j=1}^d\int_{\D^d}\langle h_j,\mathcal D_jg\rangle\\
 &=\sum_{j=1}^d\int_{\D^d}\langle\mathcal D_jh_j,g\rangle
 =\int_{\D^d}\Bigl\langle\sum_{j=1}^dh_j,g\Bigr\rangle.
\end{align*}
Applying Cauchy--Schwarz on the product of the two probability spaces and then the upper estimate in \(X^*\) bounds the absolute value of this expression by
\[
 \Bigl(\E_r\bigl\|\sum_{j=1}^dr_jh_j\bigr\|_2^2\Bigr)^{1/2}
 C_{X^*}\Bigl\|\sum_{j=1}^d\mathcal D_jg\Bigr\|_2.
\]
The last sum telescopes to \(g-\mathcal E_0g\), whose norm is at most \(2\|g\|_2\). Taking the supremum over \(\|g\|_2\leq1\) in the finite-cube duality gives
\[
 \Bigl\|\sum_{j=1}^dh_j\Bigr\|_2
 \leq2C_{X^*}
       \Bigl(\E_r\bigl\|\sum_{j=1}^dr_jh_j\bigr\|_2^2\Bigr)^{1/2}.
\]

For fixed signs \(a_1,\ldots,a_d\), the sequence \((a_jh_j)\) is again a martingale difference sequence. Applying this lower estimate to it, using the equality in distribution of \((r_ja_j)\) and \((r_j)\), and finally applying the upper estimate in \(X\), we obtain
\begin{align*}
 \Bigl\|\sum_{j=1}^da_jh_j\Bigr\|_2
 &\leq2C_{X^*}
    \Bigl(\E_r\bigl\|\sum_{j=1}^dr_ja_jh_j\bigr\|_2^2\Bigr)^{1/2}\\
 &=2C_{X^*}
    \Bigl(\E_r\bigl\|\sum_{j=1}^dr_jh_j\bigr\|_2^2\Bigr)^{1/2}
 \leq2C_{X^*}C_X\Bigl\|\sum_{j=1}^dh_j\Bigr\|_2.
\end{align*}
This bound is uniform in \(d\) and in the deterministic signs. The Walsh--Paley characterisation \cite[Theorem~4.2.13]{HNVW1} therefore implies that \(X\) is UMD. Finally, every UMD\(^+\) space has the Stein property by Bourgain's inequality, as already mentioned in Section~\ref{sec:preliminaries}. If the equivalent assertions hold, each such space is therefore UMD as well.
\end{proof}

\subsection{Stein implies the triangular contraction property}
\label{subsec:triangular-contractions}

To connect this section with the discussion of \(R\)-sectoriality in the Introduction, the final result of this section shows that if a Banach space \(X\) has the Stein property, then \(X\) has the triangular contraction property; as a consequence, every space \(E = L^p(S,\mu;X)\) with \(1\le p<\infty\) has the triangular contraction property \cite[Proposition 7.5.8]{HNVW2}.

Let \((r_i)_{i\geq0}\) be a sequence of independent symmetric signs on a probability space \(\Omega_r\), and let \((\varepsilon_j)_{j\geq0}\) be an independent Steinhaus sequence on \(\Omega_\varepsilon\). For \(1<p<\infty\), we denote by \(\Delta_{p,X}\) the least constant, if one exists, for which
\begin{equation}
\label{eq:triangular-contraction}
 \Bigl\|\sum_{j=0}^n\sum_{i=0}^j
       \varepsilon_j r_i x_{ij}\Bigr\|_{L^p(\Omega_\varepsilon\times\Omega_r;X)}
 \leq \Delta_{p,X}
 \Bigl\|\sum_{j=0}^n\sum_{i=0}^n
       \varepsilon_j r_i x_{ij}\Bigr\|_{L^p(\Omega_\varepsilon\times\Omega_r;X)}
\end{equation}
for every \(n\geq0\) and every array \((x_{ij})_{i,j=0}^n\) in \(X\). The finiteness of \(\Delta_{p,X}\) is the \emph{triangular contraction property}, or property \((\Delta)\). The usual formulation uses two Rademacher sequences and the second moment \cite[(3.6)]{KaltonWeis}. The Kahane--Khintchine inequalities applied to each sequence in turn, and the comparison of Rademacher and Steinhaus sums \cite[Section~6.1.a]{HNVW2}, show that the formulations are equivalent and that finiteness is independent of \(p\).

\begin{proposition}
\label{prop:Stein-triangular}
If \(X\) has the \(p\)-Stein property, then it has the triangular contraction property, with
\[
 \Delta_{p,X}\leq\mathfrak s_{p,X}=\mathfrak s^{\rm dyad}_{p,X}.
\]
\end{proposition}

\begin{proof}
Fix \(n\geq0\) and an array \((x_{ij})_{i,j=0}^n\) in \(X\). On the product probability space \(S:=\Omega_\varepsilon\times\Omega_r\), consider the filtration given by
\[
 \mathscr G_0:=\sigma(\varepsilon_0,\ldots,\varepsilon_n),
 \qquad
 \mathscr G_{j+1}:=\mathscr G_0\vee\sigma(r_0,\ldots,r_j)
 \quad (0\leq j\leq n).
\]
For \(0\leq j\leq n\), put
\[
 F_j:=\sum_{i=0}^n\varepsilon_j r_i x_{ij},
 \qquad
 G_j:=\sum_{i=0}^j\varepsilon_j r_i x_{ij}.
\]
The variables \(\varepsilon_j,r_0,\ldots,r_j\) are \(\mathscr G_{j+1}\)-measurable, and each remaining sign is centred and independent of \(\mathscr G_{j+1}\). Hence \(G_j=\E(F_j\mid\mathscr G_{j+1})\).

Applying the \(p\)-Stein inequality with an external Steinhaus sequence \((\eta_j)_{j=0}^n\), independent of the variables on \(S\), gives
\[
 \Bigl\|\sum_{j=0}^n\eta_j G_j\Bigr\|_{\Rad_{n+1}(L^p(S;X))}
 \leq \mathfrak s_{p,X}
 \Bigl\|\sum_{j=0}^n\eta_j F_j\Bigr\|_{\Rad_{n+1}(L^p(S;X))}.
\]
For any fixed unimodular scalars \(a_0,\ldots,a_n\), the sequence \((a_j\varepsilon_j)_{j=0}^n\) has the same joint law as \((\varepsilon_j)_{j=0}^n\) and is independent of the signs. Consequently,
\[
 \Bigl\|\sum_{j=0}^n a_j G_j\Bigr\|_{L^p(S;X)}
 =\Bigl\|\sum_{j=0}^n G_j\Bigr\|_{L^p(S;X)},
 \qquad
 \Bigl\|\sum_{j=0}^n a_j F_j\Bigr\|_{L^p(S;X)}
 =\Bigl\|\sum_{j=0}^n F_j\Bigr\|_{L^p(S;X)}.
\]
Taking \(a_j=\eta_j\) for each realisation, we see that the external \(L^2\)-norms in the Stein estimate are norms of constants. The resulting inequality is precisely \eqref{eq:triangular-contraction} with constant \(\mathfrak s_{p,X}\). Theorem~\ref{thm:full-dyadic-Stein} gives the equality with the dyadic constant.
\end{proof}
%
%

\begin{remark}
\label{rem:Stein-iterated-triangular}
There is an iterated version of \eqref{eq:triangular-contraction} in which each term has the form
\[
 \varepsilon_j r_{i_1}^{(1)}\cdots r_{i_k}^{(k)}
 x_{i_1,\ldots,i_k,j},
 \qquad 0\leq i_1,\ldots,i_k,j\leq n,
\]
where the \(k\) sign sequences are mutually independent and independent of the Steinhaus sequence. On the left, the sum retains exactly the terms with \(i_1,\ldots,i_k\leq j\). On the right, the sum runs over all indices. Both norms use the \(p\)th moment over all the sequences. It can then be shown that a Banach space has the Stein property if and only if these estimates hold with one constant for all \(n\geq0\) and \(k\geq1\).

For the forward estimate, we use the proof of Proposition~\ref{prop:Stein-triangular} with the \(\sigma\)-algebra at level \(j+1\) generated by the Steinhaus variables and the signs \(r_i^{(s)}\) with \(i\leq j\) and \(1\leq s\leq k\). For the converse on \(\D^d\), take \(n=k=d\). For each \(A\subseteq\{1,\ldots,d\}\), choose \(i_s=s\) when \(s\in A\) and \(i_s=0\) otherwise. Then
\[
 \prod_{s=1}^d r_{i_s}^{(s)}
 =\Bigl(\prod_{s=1}^d r_0^{(s)}\Bigr)
  \prod_{s\in A}(r_0^{(s)}r_s^{(s)}).
\]
The second factor is the Walsh monomial \(w_A\) evaluated at the independent symmetric signs \(r_0^{(s)}r_s^{(s)}\). Taking the norms of the two sums removes the first factor, which is common to all terms and has modulus one. The restriction \(i_1,\ldots,i_d\leq j\) selects precisely the monomials fixed by \(\mathcal E_j^{(d)}\), namely those with \(A\subseteq\{1,\ldots,j\}\).

Let \(f_1,\ldots,f_m\in L^p(\D^d;X)\) and let \(\tau(1),\ldots,\tau(m)\in\{0,\ldots,d\}\), with repetitions allowed. We take an auxiliary Steinhaus sequence \((\eta_\ell)_{\ell=1}^m\), independent of \((\varepsilon_j)_{j=0}^d\), and put
\[
 F_j:=\sum_{\ell:\,\tau(\ell)=j}\eta_\ell f_\ell,
 \qquad 0\leq j\leq d.
\]
For each fixed realisation of the auxiliary variables, we apply the estimate with one function at each level to \(F_0,\ldots,F_d\). The two sums in that estimate are
\begin{align*}
 \sum_{j=0}^d\varepsilon_j\mathcal E_j^{(d)}F_j
 =\sum_{\ell=1}^m\varepsilon_{\tau(\ell)}\eta_\ell
     \mathcal E_{\tau(\ell)}^{(d)}f_\ell,\qquad
 \sum_{j=0}^d\varepsilon_j F_j
 =\sum_{\ell=1}^m\varepsilon_{\tau(\ell)}\eta_\ell f_\ell.
\end{align*}
Conditional on \((\varepsilon_j)_{j=0}^d\), the products \(\varepsilon_{\tau(\ell)}\eta_\ell\) are independent and uniformly distributed on the unit circle. Thus taking \(p\)th powers and integrating the estimate over the auxiliary variables gives the \(p\)th-moment Stein inequality for the prescribed list of levels. The Kahane--Khintchine inequalities then give the formulation with an external second moment.
\end{remark}

\section{The Stein property and dimension-free \texorpdfstring{\(R\)}{R}-sectoriality}\label{sec:dimensionfree}

For each \(d\geq1\) we consider the following operators.
\begin{itemize}
\item The Walsh number operator \(W_d\) on \(L^p(\D^d;X)\).
\item The Ornstein--Uhlenbeck number operator \(A_d=-\Delta+x\cdot\nabla\) on \(L^p(\gamma_d;X)\).
\item The positive Laplacian \(\mathcal L_d=-\Delta\) on \(L^p(\mathbb R^d;X)\).
\item The positive Poisson operator \(\mathcal P_d=\sqrt{-\Delta}\) on \(L^p(\mathbb R^d;X)\).
\end{itemize}

In the first line, the cubes \(\D^d\) are assumed to carry uniform probability measure assigning mass \(1/2^d\) to every point. The {\em Walsh number operator} \(W_d\) on \(L^p(\D^d)\) is defined by
\[
 W_dw_A=|A|w_A,
 \qquad
 w_A(x):=\prod_{i\in A}x_i,
 \quad A\subseteq\{1,\ldots,d\}.
\]
Its extension to \(L^p(\D^d;X)\) is again denoted by \(W_d\). Here \(\gamma_d\) is standard Gaussian measure on \(\mathbb R^d\). For each family we write \((B_d,E_d(X))\) for the operator and its underlying space in dimension \(d\). In each case \(-B_d\) is the generator of the corresponding contraction semigroup. This specifies the operator domains. We write \(\Sigma_\rho=\{z\ne0:\,|\arg z|<\rho\}\).

\begin{theorem}[Stein's inequality and dimension-free \(R\)-sectoriality]
\label{thm:main}
Let \(X\) be a complex Banach space, let \(1<p<\infty\), and choose any one of the four operator families \((B_d,E_d(X))\) listed above. The following assertions are equivalent:
\begin{enumerate}[\rm(1)]
\item\label{thm:main1} \(X\) has the Stein property;
\item\label{thm:main2} the operators \(B_d\) are uniformly \(R\)-sectorial: there exists \(0<\rho<\pi/2\) such that \(\sigma(B_d)\subseteq\overline{\Sigma_\rho}\) for all \(d\ge 1\) and
\[
 \sup_{d\ge 1}\cR_{E_d(X)}
       \bigl(\{\lambda(\lambda-B_d)^{-1}:\,
              \lambda\notin\overline{\Sigma_\rho}\}\bigr)<\infty;
\]
\item\label{thm:main3} we have \((0,\infty)\subseteq \varrho(-B_d)\) and the resolvents satisfy
\[
 \sup_{d\ge 1}\cR_{E_d(X)}
       \bigl(\{(I+sB_d)^{-1}:\, s>0\}\bigr)<\infty.
\]
\end{enumerate}
If these equivalent conditions hold, then
\begin{equation}
\label{eq:resolvent-controls-Stein}
 \mathfrak s_{p,X}
 = \sup_{d\ge 1}\cR_{E_d(X)}
       \bigl(\{(I+sB_d)^{-1}:\,s>0\}\bigr).
\end{equation}
\end{theorem}

For a family \((B_d,E_d(X))\) and a parameter set \(\Lambda\), we write
\[
 C_{\mathrm{family}}(\Lambda)=
 \sup_{d\geq1}\cR_{E_d(X)}
 \bigl(\{\lambda(\lambda-B_d)^{-1}:\,\lambda\in\Lambda\}\bigr).
\]
The subscripts \(\mathrm W\), \(\mathrm{OU}\), \(\Delta\), and \(\mathrm{Poi}\) refer to the Walsh, Ornstein--Uhlenbeck, heat, and Poisson families, respectively. If a parameter in \(\Lambda\) lies outside the resolvent set for some \(d\), we set the corresponding constant equal to \(+\infty\). We omit \(\Lambda\) when it is \((-\infty,0)\), in which case \(\lambda=-1/s\) gives \(\lambda(\lambda-B_d)^{-1}=(I+sB_d)^{-1}\). The constants \(C_{\Delta,1}\) and \(C_{\mathrm{Poi},1}\) are the bounds for these families in dimension one.

The exact constants follow from the chain
\begin{equation}
\label{eq:Stein-six-bounds}
 \mathfrak s_{p,X}
 \stackrel{\rm(i)}{\leq} C_{\mathrm{Poi},1}
 \stackrel{\rm(ii)}{\leq} C_{\Delta,1}
 \stackrel{\rm(iii)}{\leq} C_\Delta
 \stackrel{\rm(iv)}{\leq} C_{\mathrm{OU}}
 \stackrel{\rm(v)}{\leq} C_{\mathrm W}
 \stackrel{\rm(vi)}{\leq} \mathfrak s_{p,X}.
\end{equation}
We prove the six inequalities on \((-\infty,0)\) in reverse order. The transfers in {\rm(v)} and {\rm(iv)} also apply to resolvents outside a sector, and we use them in that form in the completion of the proof. For the full Poisson family,
\[
 C_{\mathrm{Poi},1}\leq C_{\mathrm{Poi}}\leq C_\Delta.
\]
The first inequality holds because the supremum defining \(C_{\mathrm{Poi}}\) includes dimension one. The second follows by taking the supremum over all dimensions in the resolvent comparison proved in {\rm(ii)}. Once \eqref{eq:Stein-six-bounds} is established, both endpoints are equal to \(\mathfrak s_{p,X}\), so \(C_{\mathrm{Poi}}=\mathfrak s_{p,X}\) as well.

\subsection{Proof of {\rm (vi)}}

\begin{proof}[Proof of {\rm (vi)}]
We use the random subset representation of the Walsh semigroup and the nested conditional expectations from \cite[Section~3, (3.1)--(3.7)]{Neerven-OU}. For \(S\subseteq\{1,\ldots,d\}\) let \(M_S\) denote conditional expectation onto the \(\sigma\)-algebra generated by the coordinates with indices in \(S\). For \(u\in[0,1]^d\) and \(t\geq0\), put
\[
 S_t(u):=\{i:\,u_i\leq e^{-t}\}.
\]
Averaging a Walsh monomial over a coordinate in \(A\) gives zero, while averaging over a coordinate outside \(A\) leaves it unchanged. Thus \(M_Sw_A=\mathbf1_{\{A\subseteq S\}}w_A\).

We equip \([0,1]^d\) with product Lebesgue measure, so that the coordinate functions \(u_1,\ldots,u_d\) are independent and uniform on \([0,1]\). The inclusion \(A\subseteq S_t(u)\) holds exactly when \(u_i\leq e^{-t}\) for every \(i\in A\). Its indicator is \(\prod_{i\in A}\mathbf1_{[0,e^{-t}]}(u_i)\). Integrating one coordinate at a time gives
\[
 \int_{[0,1]^d}\mathbf1_{\{A\subseteq S_t(u)\}}\dd u
 =\prod_{i\in A}\int_0^1\mathbf1_{[0,e^{-t}]}(u_i)\dd u_i
 =e^{-t|A|}.
\]
Each coordinate in \(A\) contributes the length \(e^{-t}\), and each coordinate outside \(A\) contributes a factor of one. If \(A=\varnothing\), the inclusion holds for every \(u\) and the empty product is one. Using the formula for \(M_Sw_A\), we obtain
\[
 \int_{[0,1]^d}M_{S_t(u)}w_A\dd u
 =e^{-t|A|}w_A.
\]
Therefore, by the Walsh expansion,
\[
 e^{-tW_d}=\int_{[0,1]^d}M_{S_t(u)}\dd u.
\]
Combining this identity with the Laplace representation of the resolvent, we obtain
\[
 (I+sW_d)^{-1}
 =\int_0^\infty e^{-v}
      \int_{[0,1]^d}M_{S_{sv}(u)}\dd u\dd v,
 \qquad s>0.
\]

Fix \(s_1,\ldots,s_m>0\) and \(f_1,\ldots,f_m\in L^p(\D^d;X)\), and write \(Y_m=\Rad_m(L^p(\D^d;X))\). For each fixed \(u\) and \(v\), the sets \(S_{s_jv}(u)\) are nested after ordering the parameters \(s_j\). Their associated \(\sigma\)-algebras consequently form a finite filtration when listed in increasing order. The Stein property therefore gives
\[
 \Bigl\|\sum_{j=1}^m\varepsilon_j
                   M_{S_{s_jv}(u)}f_j\Bigr\|_{Y_m}
 \leq \mathfrak s_{p,X}
 \Bigl\|\sum_{j=1}^m\varepsilon_jf_j\Bigr\|_{Y_m}.
\]
Minkowski's inequality now implies
\begin{align*}
 \Bigl\|\sum_{j=1}^m\varepsilon_j
                   (I+s_jW_d)^{-1}f_j\Bigr\|_{Y_m}
 &\leq
 \int_0^\infty e^{-v}\int_{[0,1]^d}
 \Bigl\|\sum_{j=1}^m\varepsilon_j
                   M_{S_{s_jv}(u)}f_j\Bigr\|_{Y_m}\dd u\dd v\\
 &\leq \mathfrak s_{p,X}
 \Bigl\|\sum_{j=1}^m\varepsilon_jf_j\Bigr\|_{Y_m}.
\end{align*}
Both averaging measures have total mass one, so no additional constant occurs. This proves the assertion for the Walsh family.
\end{proof}

\subsection{Proof of {\rm (v)}}

\begin{proof}[Proof of {\rm (v)}]
We take either \(\Lambda=(-\infty,0)\) or \(\Lambda=\mathbb C\setminus\overline{\Sigma_\rho}\) with \(0<\rho<\pi/2\), and assume \(C_{\mathrm W}(\Lambda)<\infty\). The Walsh resolvents \(\lambda(\lambda-W_d)^{-1}\) have degree multipliers
\[
 m_\lambda(k)=\frac{\lambda}{\lambda-k},
 \qquad k\geq0,\quad \lambda\in\Lambda.
\]
By the Walsh-to-Hermite transfer \cite[Lemma~7.1]{Neerven-OU}, for each \(d\) the Hermite multipliers \(T_{\lambda,d}\) with symbols \(m_\lambda\) have \(R\)-bound at most \(C_{\mathrm W}(\Lambda)\) on \(L^p(\mathbb R^d,\gamma_d;X)\). The proof preserves the randomisation bound, since the convergence holds for each fixed choice of unimodular coefficients and the uniform moment bounds justify integration over the coefficients.

To identify \(T_{\lambda,d}\), we use that \(I+A_d\) is invertible because \(-A_d\) generates a contraction semigroup. The bounded operators
\[
 (\lambda-A_d)(I+A_d)^{-1}
 \qquad\text{and}\qquad
 \frac{1+\lambda}{\lambda}T_{\lambda,d}-I
\]
act on Hermite terms of degree \(k\) by the reciprocal scalars \((\lambda-k)/(1+k)\) and \((1+k)/(\lambda-k)\). By density of the finite Hermite expansions, these operators are inverses on \(L^p(\mathbb R^d,\gamma_d;X)\). Since \(I+A_d\) maps \(\Dom(A_d)\) bijectively onto this space, it follows that \(\lambda\in\varrho(A_d)\). The normalised resolvent \(\lambda(\lambda-A_d)^{-1}\) agrees with \(T_{\lambda,d}\) on finite Hermite expansions and hence on the whole space.

Taking the supremum over \(d\) gives \(C_{\mathrm{OU}}(\Lambda)\leq C_{\mathrm W}(\Lambda)\). For \(\Lambda=(-\infty,0)\) this is {\rm(v)}.
\end{proof}

\subsection{Proof of {\rm (iv)}}

\begin{proof}[Proof of {\rm (iv)}] Fix \(d\ge 1\).
Let \(\Lambda\) be either of the sets in {\rm(v)}. We may assume that \(C_{\mathrm{OU}}(\Lambda)<\infty\). The Laplacian \(\mathcal L_d\) is sectorial of angle zero on \(L^p(\mathbb R^d;X)\), so \(\Lambda\subseteq\varrho(\mathcal L_d)\) \cite[Example~10.1.5]{HNVW2}.

We use the change of density and dilation from the one-dimensional \(H^\infty\)-calculus comparison in \cite[proof of Theorem~9.1, (9.2)--(9.4)]{Neerven-OU}. Their combined action is described by the family of operators \(V_a : L^p(\mathbb R^d;X) \to L^p(\mathbb R, \gamma_d;X)\), \(a\ge 1\), given by
\[
 V_ag(x):=(2\pi)^{d/(2p)}a^{d/p}e^{|x|^2/(2p)}g(ax).
\]
The exponential weight cancels the Gaussian density in the norm, and the substitution \(y=ax\) gives
\[
 \|V_ag\|_{L^p(\mathbb R^d,\gamma_d;X)}^p
 =a^d\int_{\mathbb R^d}\|g(ax)\|_X^p\dd x
 =\|g\|_{L^p(\mathbb R^d;X)}^p.
\]
Reversing the dilation and dividing by the weight gives the inverse, so \(V_a\) is a surjective isometry which acts bijectively on \(\mathscr D_d:= C_{\rm c}^\infty(\mathbb R^d)\otimes X\). We define
\[
 L_{d,a}:=a^{-2}V_a^{-1}A_dV_a,
 \qquad \Dom(L_{d,a})=V_a^{-1}\Dom(A_d).
\]
We have \(\mathscr D_d \subseteq \Dom(A_d)\),
and on \(\mathscr D_d\) the operator \(A_d\) is given as \( -\Delta+x\cdot\nabla\); in particular, \(A_d\) maps \(\mathscr D_d\) into itself. It follows that \(\mathscr D_d\subseteq \Dom(L_{d,a})\).

For the weight \(w(y)=e^{|y|^2/(2p)}\), we have \(\nabla w(y)=p^{-1}yw(y)\) and \(\Delta w(y)=(d/p+|y|^2/p^2)w(y)\). Applying the product rule to \(A_d(wg(a\,\cdot))\) and evaluating at \(y=x/a\) in the formula for \(V_a^{-1}\) gives the explicit expression
\[
 L_{d,a}g
 =-\Delta g+\frac{1-2/p}{a^2}x\cdot\nabla g
   +\frac{p-1}{p^2a^4}|x|^2g-\frac{d}{pa^2}g, \qquad g\in \mathscr D_d.
\]
For fixed \(g\in\mathscr D_d\), the three lower-order terms tend to zero in \(L^p(\mathbb R^d;X)\) as \(a\to\infty\). Thus, for \(g\in \mathscr D_d\), we have \(L_{d,a}g\to\mathcal L_dg\) as \(a\to\infty\).

For \(\lambda\in\Lambda\), we have \(a^2\lambda\in\Lambda\) by invariance under positive dilations.
Since \(\Lambda\subseteq \varrho(A_d)\) it follows that
\(a^2\lambda\in\varrho(A_d)\).
The definition of \(L_{d,a}\) therefore gives
\[
 R_a(\lambda):=\lambda(\lambda-L_{d,a})^{-1}
 =V_a^{-1}\bigl[a^2\lambda(a^2\lambda-A_d)^{-1}\bigr]V_a.
\]
As an immediate consequence, the \(R\)-boundedness of the operators \(a^2\lambda(a^2\lambda-A_d)^{-1}\) implies the \(R\)-boundedness of the operators \(R_a(\lambda)\), with \(R\)-bound
\begin{align}\label{eq:Ra-Rbdd}
 \cR_{L^p(\mathbb R^d;X)}\bigl(\{R_a(\lambda):\,\lambda\in\Lambda\}\bigr)
 \leq C_{\mathrm{OU}}(\Lambda)
 \qquad (a\geq1).
\end{align}

Fix \(\lambda\in\Lambda\) and put \(R(\lambda)=\lambda(\lambda-\mathcal L_d)^{-1}\). If \(g\in\mathscr D_d\) and \(f=(\lambda-\mathcal L_d)g\), then
\[
 R_a(\lambda)f-R(\lambda)f
 = R_a(\lambda)((\lambda-L_{d,a} + L_{d,a} - \mathcal L_d)g-\lambda g
 =R_a(\lambda)(L_{d,a}-\mathcal L_d)g.
\]
The right-hand side tends to zero in \(L^p(\mathbb R^d;X)\), since \(L_{d,a}g\to\mathcal L_dg\) and \(\|R_a(\lambda)\|\leq C_{\mathrm{OU}}(\Lambda)\).
This proves the convergence
\[
 R_a(\lambda)f-R(\lambda)f \to 0, \qquad f \in  (\lambda-\mathcal L_d)\mathscr D_d.
\]
Since \(C_{\rm c}^\infty(\mathbb R^d)\) is a core for \(\lambda+\Delta\), \((\lambda+\Delta)C_{\rm c}^\infty(\mathbb R^d)\) is dense in \(L^p(\mathbb R^d)\). Tensoring with \(X\) shows that \((\lambda-\mathcal L_d)\mathscr D_d\) is dense in \(L^p(\mathbb R^d;X)\). Since \(R(\lambda)\) is bounded and the bound for \(R_a(\lambda)\) is uniform in \(a\),
we conclude that
\[
 R_a(\lambda)f\to R(\lambda)f, \qquad f\in L^p(\mathbb R^d;X).
\]

Finally, let \(m\geq1\), \(\lambda_1,\ldots,\lambda_m\in\Lambda\) and \(f_1,\ldots,f_m\in L^p(\mathbb R^d;X)\). For every \(a\geq1\), from \eqref{eq:Ra-Rbdd} we have
\[
 \Bigl\|\sum_{\ell=1}^m\varepsilon_\ell R_a(\lambda_\ell)f_\ell
       \Bigr\|_{\Rad_m(L^p(\mathbb R^d;X))}
 \leq C_{\mathrm{OU}}(\Lambda)
 \Bigl\|\sum_{\ell=1}^m\varepsilon_\ell f_\ell\Bigr\|_{\Rad_m(L^p(\mathbb R^d;X))}.
\]
For each \(\ell\), we have \(R_a(\lambda_\ell)f_\ell\to R(\lambda_\ell)f_\ell\) in \(L^p(\mathbb R^d;X)\) as \(a\to\infty\). The finite sum on the left therefore converges in \(\Rad_m(L^p(\mathbb R^d;X))\). Letting \(a\to\infty\) in the preceding inequality gives
\[
 \Bigl\|\sum_{\ell=1}^m\varepsilon_\ell R(\lambda_\ell)f_\ell
       \Bigr\|_{\Rad_m(L^p(\mathbb R^d;X))}
 \leq C_{\mathrm{OU}}(\Lambda)
 \Bigl\|\sum_{\ell=1}^m\varepsilon_\ell f_\ell
       \Bigr\|_{\Rad_m(L^p(\mathbb R^d;X))}.
\]
Since \(m\), the parameters \(\lambda_\ell\in\Lambda\) and the functions \(f_\ell\in L^p(\mathbb R^d;X)\) were arbitrary, this proves
\[
 \cR_{L^p(\mathbb R^d;X)}
 \bigl(\{\lambda(\lambda-\mathcal L_d)^{-1}:\,\lambda\in\Lambda\}\bigr)
 \leq C_{\mathrm{OU}}(\Lambda).
\]
Taking the supremum over \(d\) gives \(C_\Delta(\Lambda)\leq C_{\mathrm{OU}}(\Lambda)\). In particular, \(\Lambda=(-\infty,0)\) gives {\rm(iv)}.
\end{proof}

\subsection{Proof of {\rm (iii)}}

\begin{proof}[Proof of {\rm (iii)}] This inequality is immediate from the definitions.
\end{proof}

\subsection{Proof of {\rm (ii)}}

\begin{proof}[Proof of {\rm (ii)}]
Fix \(d\geq1\). Since \(\mathcal L_d\) is sectorial of angle zero, its square root \(\mathcal P_d\) is sectorial of angle zero as well by \cite[Theorem~10.3]{Komatsu-Fractional}.

For \(t>0\), the square-root case of Kato's resolvent formula \cite[Theorem~2, (3)]{Kato-Fractional-Note}, also given in \cite[(1.7)]{Kato-Fractional}, reads
\begin{equation*}
 t(t+\mathcal P_d)^{-1}
 =\frac1\pi\int_0^\infty
   \frac{t\,\tau^{-1/2}}{t^2+\tau}\,
   \tau(\tau+\mathcal L_d)^{-1}\dd\tau.
\end{equation*}
Taking \(t=1/s\) and substituting \(\tau=u^2/s^2\) gives
\[
 (I+s\mathcal P_d)^{-1}
 =\int_0^\infty (I+(s/u)^2\mathcal L_d)^{-1}
       \frac{2\,\dd u}{\pi(1+u^2)}.
\]
The measure on the right is a probability measure independent of \(s\). In a finite random sum we can therefore use the same variable \(u\) for every resolvent parameter. Minkowski's inequality and the heat resolvent \(R\)-bound at each \(u\) give the estimate with the same constant, as in the integral means theorem \cite[Theorem~8.5.2]{HNVW2}. Thus
\begin{equation*}
 \cR_{L^p(\mathbb R^d;X)}
   \bigl(\{(I+s\mathcal P_d)^{-1}:\,s>0\}\bigr)
 \leq
 \cR_{L^p(\mathbb R^d;X)}
   \bigl(\{(I+s\mathcal L_d)^{-1}:\,s>0\}\bigr).
\end{equation*}
The case \(d=1\) gives {\rm(ii)}, and taking the supremum over \(d\) gives \(C_{\mathrm{Poi}}\leq C_\Delta\).
\end{proof}

\subsection{Proof of {\rm (i)}}
\label{sec:onedim-recovery}

We first show that the \(R\)-bound for the Poisson resolvents on the real line also holds on the circle. We then compose Fourier characters on a product torus with a map from the circle. For each \(j\), we choose the resolvent parameters so that the multipliers on these characters converge to those of conditional expectation onto the first \(j\) coordinates.

\subsubsection{Periodisation without loss of the constant}

Let \(\mathbb T=\mathbb R/(2\pi\mathbb Z)\), let \(m_n\) denote normalised Haar measure on \(\mathbb T^n\), and write \(\widetilde F\) for the periodic extension of a function \(F\) on \(\mathbb T\). The following convolution form of de Leeuw's restriction principle \cite{deLeeuw} preserves the constant in the finite random sums defining the \(R\)-bound.

\begin{lemma}[Periodisation of convolution families]
\label{lem:onedim-periodisation}
Let \((k_\alpha)_{\alpha\in I}\) be a family in \(L^1(\mathbb R)\), and put \(K_\alpha f=k_\alpha*f\) on \(L^p(\mathbb R;X)\). The operators
\[
 K_\alpha^{\mathbb T}F(\theta)=\int_{\mathbb R} k_\alpha(u)\widetilde F(\theta-u)\dd u,
\]
initially defined for continuous \(F\), extend to \(L^p(\mathbb T;X)\) and satisfy
\begin{equation}
\label{eq:onedim-periodisation-bound}
 \cR_{L^p(\mathbb T;X)}\bigl(\{K_\alpha^{\mathbb T}:\alpha\in I\}\bigr)
 \leq \cR_{L^p(\mathbb R;X)}\bigl(\{K_\alpha:\alpha\in I\}\bigr),
\end{equation}
where the right-hand side is allowed to be infinite.
Moreover, \(K_\alpha^{\mathbb T}(e^{ik\cdot}x)=\widehat k_\alpha(k)e^{ik\cdot}x\) for all \(k\in\mathbb Z\) and \(x\in X\), where \(\widehat k_\alpha(\xi)=\int_{\mathbb R}e^{-iu\xi}k_\alpha(u)\dd u\).
\end{lemma}

\begin{proof}
Minkowski's inequality and translation invariance of \(m_1\) give \(\|K_\alpha^{\mathbb T}F\|_{p}\leq\|k_\alpha\|_1\|F\|_{p}\), so the formula extends by density. Substitution of \(F(\theta)=e^{ik\theta}x\) gives the stated multiplier.

For an integer \(N\geq1\), put \(Q_N=(-\pi N,\pi N]\), \(I_NF=|Q_N|^{-1/p}\mathbf1_{Q_N}\widetilde F\) and \(P_Nf=\mathbf1_{Q_N}f\). The interval \(Q_N\) consists of \(N\) periods of length \(2\pi\), so
\[
 \|I_NF\|_{L^p(\mathbb R;X)}^p
 =\frac{1}{|Q_N|}\int_{Q_N}\|\widetilde F(u)\|_X^p\dd u
 =\int_{\mathbb T}\|F(\theta)\|_X^p\dd m_1(\theta).
\]
Thus \(I_N\) is an isometry, and \(P_N\) is contractive.

Let \(F\) be continuous. For \(u\in Q_N\), the convolution defining \(P_NK_\alpha I_NF(u)\) includes only those \(y\in\mathbb R\) for which \(u-y\in Q_N\), whereas the integral defining \(I_NK_\alpha^{\mathbb T}F(u)\) is over all of \(\mathbb R\). Both expressions vanish for \(u\notin Q_N\). Since \(u-y\in Q_N\) is equivalent to \(u\in Q_N+y\), their difference is
\begin{align*}
 &(P_NK_\alpha I_NF-I_NK_\alpha^{\mathbb T}F)(u)\\
 &\qquad=-|Q_N|^{-1/p}\int_{\mathbb R} k_\alpha(y)
   \mathbf1_{Q_N\setminus(Q_N+y)}(u)\widetilde F(u-y)\dd y.
\end{align*}
For fixed \(y\), the \(L^p(\mathbb R;X)\)-norm of the product of the indicator and \(\widetilde F(\,\cdot-y)\) is at most \(\|F\|_\infty |Q_N\setminus(Q_N+y)|^{1/p}\). Minkowski's integral inequality therefore gives
\begin{align*}
 &\|P_NK_\alpha I_NF-I_NK_\alpha^{\mathbb T}F\|_{L^p(\mathbb R;X)}\\
 &\qquad\leq \|F\|_\infty
 \int_{\mathbb R}|k_\alpha(y)|
 \left(\frac{|Q_N\setminus(Q_N+y)|}{|Q_N|}\right)^{1/p}\dd y.
\end{align*}
When \(|y|<2\pi N\), the two intervals overlap in an interval of length \(2\pi N-|y|\). When \(|y|\geq2\pi N\), their intersection has measure zero. In either case,
\[
 \frac{|Q_N\setminus(Q_N+y)|}{|Q_N|}
 =\min\left\{1,\frac{|y|}{2\pi N}\right\}.
\]
This ratio tends to zero for each fixed \(y\) and is bounded by one. Since \(k_\alpha\in L^1(\mathbb R)\), dominated convergence shows that the norm of the difference tends to zero.

Fix a finite subset \(J\subseteq I\) and put
\[
 C_J=\cR_{L^p(\mathbb R;X)}\bigl(\{K_\alpha:\alpha\in J\}\bigr).
\]
Since \(J\) is finite, \(C_J<\infty\).

Now fix \(m\geq1\), functions \(F_1,\ldots,F_m\in C(\mathbb T;X)\),
and indices \(\alpha_1,\ldots,\alpha_m\in J\), with repetitions allowed.
The contractivity of \(P_N\), the definition of \(C_J\), and the
isometry of \(I_N\) give
\[
 \Bigl\|\sum_{\ell=1}^m\varepsilon_\ell
       P_NK_{\alpha_\ell}I_NF_\ell\Bigr\|_{\Rad_m(L^p(\mathbb R;X))}
 \leq C_J
 \Bigl\|\sum_{\ell=1}^m\varepsilon_\ell
       F_\ell\Bigr\|_{\Rad_m(L^p(\mathbb T;X))}.
\]
The preceding \(L^p\)-approximation and the triangle inequality show that the left-hand norm differs by a quantity tending to zero from the norm of \(\sum_\ell\varepsilon_\ell I_NK_{\alpha_\ell}^{\mathbb T}F_\ell\). By the isometry of \(I_N\), the latter equals the randomised norm of \(\sum_\ell\varepsilon_\ell K_{\alpha_\ell}^{\mathbb T}F_\ell\) on the circle. Letting \(N\to\infty\) gives
\[
 \Bigl\|\sum_{\ell=1}^m\varepsilon_\ell
       K_{\alpha_\ell}^{\mathbb T}F_\ell
       \Bigr\|_{\Rad_m(L^p(\mathbb T;X))}
 \leq C_J
 \Bigl\|\sum_{\ell=1}^m\varepsilon_\ell
       F_\ell\Bigr\|_{\Rad_m(L^p(\mathbb T;X))}.
\]
By density of \(C(\mathbb T;X)\) and boundedness of the
\(K_{\alpha_\ell}^{\mathbb T}\), the estimate extends to arbitrary
 \(F_\ell\in L^p(\mathbb T;X)\).

We have thus shown that the \(R\)-bound of \(\{K_\alpha^{\mathbb T}:\alpha\in J\}\) is at most \(C_J\). Every finite list of indices from \(I\) is contained in a finite subset \(J\), so taking the supremum over such \(J\) gives \eqref{eq:onedim-periodisation-bound}.
\end{proof}

Let \((p_t)_{t>0}\) be the Poisson kernel on \(\mathbb R\),
\[ p_t(x) = \frac1\pi \frac{t}{t^2+x^2},\qquad  t>0, \ x\in\mathbb R.\] With our Fourier transform convention, \(p_t\) is non-negative, has integral one and satisfies \(\widehat p_t(\xi)=e^{-t|\xi|}\) \cite[(13.26)--(13.27)]{Neerven-FA}. By the Laplace formula for the resolvent \cite[Proposition~13.8]{Neerven-FA}, the operator \((I+s\mathcal P_1)^{-1}\) is convolution with
\begin{equation}
\label{eq:onedim-Poisson-resolvent-kernel}
 q_s=\int_0^\infty e^{-u}p_{su}\dd u,
 \qquad
 \widehat q_s(\xi)=\int_0^\infty e^{-u}e^{-su|\xi|}\dd u
                 =\frac{1}{1+s|\xi|}.
\end{equation}
The integral defining \(q_s\) converges in \(L^1(\mathbb R)\), and \(q_s\) is non-negative with integral one.

Let \(\mathscr Q_s\) denote the periodisation of convolution with \(q_s\). Its multiplier is
\begin{equation}
\label{eq:onedim-circle-multiplier}
 \mathscr Q_s(e^{ik\theta}x)=\frac{e^{ik\theta}x}{1+s|k|}
 \qquad(k\in\mathbb Z,\ x\in X).
\end{equation}
This is the resolvent of the circle Poisson generator, since the circle Poisson semigroup has multiplier \(e^{-t|k|}\). Lemma~\ref{lem:onedim-periodisation} gives
\begin{equation}
\label{eq:onedim-real-to-circle}
 \cR_{L^p(\mathbb T;X)}\bigl(\{\mathscr Q_s:s>0\}\bigr)
 \leq \cR_{L^p(\mathbb R;X)}\bigl(\{(I+s\mathcal P_1)^{-1}:s>0\}\bigr).
\end{equation}

\subsubsection{\texorpdfstring{Conditional expectations on the product torus}{Conditional expectations on the product torus}}

We use a mapping from \(\mathbb T\) to \(\mathbb T^n\) whose \(i\)th coordinate is multiplication by \(b^i\). This separation of frequencies is used in Bourgain's argument for the Hilbert transform \cite[proof of Lemma~1]{Bourgain83} and in Hyt\"onen's argument for the Poisson square function \cite[Section~3]{Hytonen-LPS}.

For integers \(n\geq1\) and \(b\geq2\), define \(\beta_b:\mathbb T\to\mathbb T^n\) and \(J_b:C(\mathbb T^n;X)\to C(\mathbb T;X)\) by
\[
 \beta_b(\theta)=(b\theta,b^2\theta,\ldots,b^n\theta),
 \qquad J_bF=F\circ\beta_b.
\]
For \(k=(k_1,\ldots,k_n)\in\mathbb Z^n\), set \(\nu_b(k)=\sum_{i=1}^n k_i b^i\). Composition with \(\beta_b\) sends the character \(\eta\mapsto e^{ik\cdot\eta}\) of \(\mathbb T^n\) to the character \(\theta\mapsto e^{i\nu_b(k)\theta}\) of \(\mathbb T\). If \(0\ne k\in \mathbb Z^n\) and \(j\) is the largest index with \(k_j\ne0\), then
\begin{equation}
\label{eq:onedim-leading-frequency}
 b^{-j}\nu_b(k)=k_j+\sum_{i<j}k_i b^{i-j}\to k_j \qquad (b\to\infty).
\end{equation}

\begin{lemma}[Convergence of the randomised norms]
\label{lem:onedim-norm-convergence}
For every finite sequence \(F_1,\ldots,F_m\in C(\mathbb T^n;X)\),
\begin{equation}
\label{eq:onedim-Rad-norm-limit}
 \Bigl\|\sum_{\ell=1}^m\varepsilon_\ell J_bF_\ell\Bigr\|_{\Rad_m(L^p(\mathbb T;X))}
 \to
 \Bigl\|\sum_{\ell=1}^m\varepsilon_\ell F_\ell\Bigr\|_{\Rad_m(L^p(\mathbb T^n;X))}
 \qquad(b\to\infty).
\end{equation}
\end{lemma}

\begin{proof}
We first show that for every scalar-valued \(g\in C(\mathbb T^n)\),
\begin{align}\label{eq:conv-b}
 \int_{\mathbb T}g(\beta_b(\theta))\dd m_1(\theta)
 \to
 \int_{\mathbb T^n}g(\eta)\dd m_n(\eta).
\end{align}
For a character \(g(\eta)=e^{ik\cdot\eta}\), the integral on the left is
\[
 \int_{\mathbb T}e^{i\nu_b(k)\theta}\dd m_1(\theta)
 =
 \begin{cases}
  1, & \nu_b(k)=0,\\
  0, & \nu_b(k)\ne0.
 \end{cases}
\]
If \(k=0\), both integrals are one. If \(k\ne0\), then \eqref{eq:onedim-leading-frequency} shows that \(\nu_b(k)\ne0\) for all sufficiently large \(b\). The integral on the left is therefore eventually zero, which is also the integral of the character against \(m_n\). A trigonometric polynomial involves only finitely many characters, so for each such polynomial the two integrals agree for all sufficiently large \(b\). The Stone--Weierstrass theorem gives uniform density of these polynomials in \(C(\mathbb T^n)\). Since each integral in \eqref{eq:conv-b} is bounded by \(\|g\|_\infty\), the convergence extends to all \(g\in C(\mathbb T^n)\).

For fixed unimodular scalars \(a_1,\ldots,a_m\), the function
\[
 g_a(\eta)=\Bigl\|\sum_{\ell=1}^m a_\ell F_\ell(\eta)\Bigr\|_X^p
\]
is continuous. Applying \eqref{eq:conv-b} to \(g_a\) gives
\begin{align*}
 \Bigl\|\sum_{\ell=1}^m a_\ell J_bF_\ell\Bigr\|_{L^p(\mathbb T;X)}^p
 &=\int_{\mathbb T}g_a(\beta_b(\theta))\dd m_1(\theta)\\
 &\to\int_{\mathbb T^n}g_a(\eta)\dd m_n(\eta)
 =\Bigl\|\sum_{\ell=1}^m a_\ell F_\ell\Bigr\|_{L^p(\mathbb T^n;X)}^p.
\end{align*}
Both \(L^p\)-norms are bounded by \(\sum_{\ell=1}^m\|F_\ell\|_\infty\), independently of \(b\) and the scalars \(a_\ell\). Taking \(a_\ell=\varepsilon_\ell(\omega)\) and applying dominated convergence to the squares of the norms gives
\[
 \E \Bigl\|\sum_{\ell=1}^m \varepsilon_\ell J_bF_\ell\Bigr\|_{L^p(\mathbb T;X)}^2 \to
 \E \Bigl\|\sum_{\ell=1}^m \varepsilon_\ell F_\ell\Bigr\|_{L^p(\mathbb T^n;X)}^2 \qquad (b\to\infty).
\]
Taking square roots gives \eqref{eq:onedim-Rad-norm-limit}.
\end{proof}

For \(0\leq j\leq n\), let \(E_j\) be the conditional expectation on \(\mathbb T^n\) onto the first \(j\) coordinates. Integrating the remaining coordinates against product Haar measure gives
\begin{equation}
\label{eq:onedim-coordinate-character}
 E_j(e^{ik\cdot\theta}x)=\mathbf1_{\{k_{j+1}=\cdots=k_n=0\}}e^{ik\cdot\theta}x .
\end{equation}
For \(j=0\), this retains only the constant character, and for \(j=n\) the condition is vacuous.

\begin{proposition}[Recovery from the circle Poisson resolvents]
\label{prop:onedim-circle-recovery}
For every complex Banach space \(X\) and every \(1<p<\infty\),
\begin{equation}
\label{eq:onedim-circle-controls-Stein}
 \mathfrak s_{p,X}\leq \cR_{L^p(\mathbb T;X)}\bigl(\{\mathscr Q_s:s>0\}\bigr).
\end{equation}
\end{proposition}

\begin{proof}
Let \(C\) denote the \(R\)-bound on the right-hand side of \eqref{eq:onedim-circle-controls-Stein}. We may assume that \(C<\infty\). Fix \(n\geq1\). We first show that the conditional expectations \(E_0,\ldots,E_n\) on \(\mathbb T^n\) satisfy Stein's inequality with constant \(C\).

For \(0\leq i\leq n\), put \(s_{i,b}=1/b^{i+1/2}\). If \(0\ne k\in\mathbb Z^n\) and \(j\) is the largest index for which \(k_j\ne0\), then
\[
 s_{i,b}|\nu_b(k)|
 =b^{j-i-1/2}\bigl|b^{-j}\nu_b(k)\bigr|.
\]
By \eqref{eq:onedim-leading-frequency}, \(\bigl|b^{-j}\nu_b(k)\bigr|\to|k_j|>0\) as \(b\to\infty\). The product therefore tends to zero if \(j\leq i\) and to infinity if \(j>i\). For \(k=0\) the product is zero for every \(b\). Thus for all \(k\in\mathbb Z^n\),
\begin{equation}
\label{eq:onedim-multiplier-cutoff}
 \frac{1}{1+s_{i,b}|\nu_b(k)|}
 \to \mathbf1_{\{k_{i+1}=\cdots=k_n=0\}}
 \qquad(b\to\infty).
\end{equation}

Let \(F(\eta)=\sum_{k\in K}e^{ik\cdot\eta}x_k\) be an \(X\)-valued trigonometric polynomial on \(\mathbb T^n\), with \(K\subseteq\mathbb Z^n\) finite.
Applying \eqref{eq:onedim-circle-multiplier} to the identity
\(J_bF(\theta)=\sum_{k\in K}e^{i\nu_b(k)\theta}x_k\)
gives
\[
 (\mathscr Q_{s_{i,b}}J_bF)(\theta)
 =\sum_{k\in K}
   \frac{1}{1+s_{i,b}|\nu_b(k)|}\,e^{i\nu_b(k)\theta}x_k,
\]
and applying \(J_b\) to the expansion of \(E_iF\) obtained from
\eqref{eq:onedim-coordinate-character} gives
\[
 (J_bE_iF)(\theta)
 =\sum_{k\in K}\mathbf1_{\{k_{i+1}=\cdots=k_n=0\}}
   e^{i\nu_b(k)\theta}x_k.
\]
It follows that
\begin{align}\label{eq:QisbJb0}
 \|\mathscr Q_{s_{i,b}}J_bF-J_bE_iF\|_\infty
 \leq \sum_{k\in K}
 \left|\frac{1}{1+s_{i,b}|\nu_b(k)|}
       -\mathbf1_{\{k_{i+1}=\cdots=k_n=0\}}\right|\|x_k\|_X
 \to0.
\end{align}
Each summand tends to zero by \eqref{eq:onedim-multiplier-cutoff}.

Now fix \(m\geq1\), trigonometric polynomials \(F_1,\ldots,F_m\) on \(\mathbb T^n\), and indices \(i_1,\ldots,i_m\in\{0,\ldots,n\}\), with repetitions allowed. By the definition of \(C\),
\[
 \Bigl\|\sum_{\ell=1}^m\varepsilon_\ell
       \mathscr Q_{s_{i_\ell,b}}J_bF_\ell
       \Bigr\|_{\Rad_m(L^p(\mathbb T;X))}
 \leq C
 \Bigl\|\sum_{\ell=1}^m\varepsilon_\ell J_bF_\ell
       \Bigr\|_{\Rad_m(L^p(\mathbb T;X))}.
\]
By \eqref{eq:QisbJb0} and the triangle inequality for the finite random sum, replacing \(\mathscr Q_{s_{i_\ell,b}}J_bF_\ell\) by \(J_bE_{i_\ell}F_\ell\) changes the norm on the left by a quantity tending to zero. Lemma~\ref{lem:onedim-norm-convergence} applies to the two lists of fixed continuous functions \(F_\ell\) and \(E_{i_\ell}F_\ell\). Letting \(b\to\infty\) gives
\begin{equation}
\label{eq:onedim-product-Stein}
 \Bigl\|\sum_{\ell=1}^m\varepsilon_\ell E_{i_\ell}F_\ell
       \Bigr\|_{\Rad_m(L^p(\mathbb T^n;X))}
 \leq C
 \Bigl\|\sum_{\ell=1}^m\varepsilon_\ell F_\ell
       \Bigr\|_{\Rad_m(L^p(\mathbb T^n;X))}.
\end{equation}
By the density of trigonometric polynomials and the contractivity of the \(E_{i_\ell}\), the estimate extends to arbitrary \(F_\ell\in L^p(\mathbb T^n;X)\).

To pass to the dyadic cube, put \(r(\theta)=\operatorname{sgn}(\cos\theta)\), with arbitrary sign values at the two zeros, as in \cite[proof of Lemma~2]{Bourgain83}. Under \(m_n\) the signs \(r(\theta_1),\ldots,r(\theta_n)\) are independent and symmetric, so each point of \(\D^n\) has probability \(2^{-n}\). Hence the mapping \[Jh(\theta):=h(r(\theta_1),\ldots,r(\theta_n))\] is an isometry from \(L^p(\D^n;X)\) into \(L^p(\mathbb T^n;X)\). Here \(\mathcal E_i^{(n)}\) denotes conditional expectation on the uniform cube \(\D^n\) with respect to the \(\sigma\)-algebra generated by its first \(i\) coordinates. For \(0\leq i<n\), integration in the last \(n-i\) coordinates gives
\[
 (E_iJh)(\theta)
 =2^{-(n-i)}\sum_{\sigma\in\D^{n-i}}
       h(r(\theta_1),\ldots,r(\theta_i),\sigma)
 =(J\mathcal E_i^{(n)}h)(\theta).
\]
For \(i=n\), both sides of this intertwining identity equal \(Jh\). For \(h_1,\ldots,h_m\in L^p(\D^n;X)\) and fixed unimodular scalars \(a_1,\ldots,a_m\), linearity and the isometry of \(J\) give
\begin{align*}
 \Bigl\|\sum_{\ell=1}^m a_\ell Jh_\ell\Bigr\|_{L^p(\mathbb T^n;X)}
 = \Bigl\|J\Bigl(\sum_{\ell=1}^m a_\ell h_\ell\Bigr)\Bigr\|_{L^p(\mathbb T^n;X)}
 = \Bigl\|\sum_{\ell=1}^m a_\ell h_\ell\Bigr\|_{L^p(\D^n;X)}.
\end{align*}
The equality holds in particular for \(a_\ell = \varepsilon_\ell(\omega)\). Taking square expectations with respect to \(\omega\) gives
\[
 \Bigl\|\sum_{\ell=1}^m\varepsilon_\ell Jh_\ell\Bigr\|_{\Rad_m(L^p(\mathbb T^n;X))}
 =
 \Bigl\|\sum_{\ell=1}^m\varepsilon_\ell h_\ell\Bigr\|_{\Rad_m(L^p(\D^n;X))}.
\]
Hence, for arbitrary \(h_1,\ldots,h_m\in L^p(\D^n;X)\),
\begin{align*}
 \Bigl\|\sum_{\ell=1}^m\varepsilon_\ell
       \mathcal E_{i_\ell}^{(n)}h_\ell\Bigr\|_{\Rad_m(L^p(\D^n;X))}
 &=\Bigl\|\sum_{\ell=1}^m\varepsilon_\ell
       J\mathcal E_{i_\ell}^{(n)}h_\ell\Bigr\|_{\Rad_m(L^p(\mathbb T^n;X))}\\
 &=\Bigl\|\sum_{\ell=1}^m\varepsilon_\ell
       E_{i_\ell}Jh_\ell\Bigr\|_{\Rad_m(L^p(\mathbb T^n;X))}\\
 &\leq C\Bigl\|\sum_{\ell=1}^m\varepsilon_\ell
       Jh_\ell\Bigr\|_{\Rad_m(L^p(\mathbb T^n;X))}\\
 &=C\Bigl\|\sum_{\ell=1}^m\varepsilon_\ell
       h_\ell\Bigr\|_{\Rad_m(L^p(\D^n;X))}.
\end{align*}
The first and last equalities use the preservation of the randomised
norm by \(J\). The second equality uses \(J\mathcal E_i^{(n)}=E_iJ\),
and the inequality follows from \eqref{eq:onedim-product-Stein} with
\(F_\ell=Jh_\ell\).

Since \(n\), \(m\) and the indices \(i_\ell\) were arbitrary, this gives \(\mathfrak s_{p,X}^{\rm dyad}\leq C\). Theorem~\ref{thm:full-dyadic-Stein} identifies the dyadic and full Stein constants and completes the proof.
\end{proof}

Combining Proposition~\ref{prop:onedim-circle-recovery} with \eqref{eq:onedim-real-to-circle} proves {\rm(i)}.

\begin{proof}[Completion of the proof of Theorem~\ref{thm:main}]
If \(X\) has the Stein property, the products of lifted Walsh projections in \cite[(3.6)]{Neerven-OU} are uniformly bounded. This follows from Stein's inequality because each product of the \(M_S\) in {\rm(vi)} conditions on the intersection of the coordinate sets, and these intersections are nested in the order of the time parameters. Also, \(X\) is \(K\)-convex by \cite[Lemma~4]{Eskenazis}, so the lifted spaces have a common non-trivial type exponent and uniform type constants. The uniform bounds for the products of projections and the uniform type estimates for the lifted spaces allow us to apply the argument of \cite[Proposition~4.2 and Theorems~5.1 and~6.1]{Neerven-OU}. By the equivalence of Rademacher and Steinhaus randomisation, it gives \(C_{\mathrm W}(\Lambda)<\infty\) for \(\Lambda=\mathbb C\setminus\overline{\Sigma_\rho}\) and some \(0<\rho<\pi/2\). Parts {\rm(v)} and {\rm(iv)} give \(C_\Delta(\Lambda)\leq C_{\mathrm{OU}}(\Lambda)\leq C_{\mathrm W}(\Lambda)\).

For the Poisson family, the bound \(C_\Delta(\Lambda)\) and \cite[Proposition~10.3.2 and Example~10.3.5]{HNVW2} give a dimension-free \(R\)-bound for \(\{|\xi|\mathcal P_d(\xi^2+\mathcal L_d)^{-1}:\xi\in\mathbb R\setminus\{0\}\}\). The factorisation of \(-(\xi^2+\mathcal L_d)\) using \(\mathcal P_d^2=\mathcal L_d\) gives
\[
 i\xi(i\xi-\mathcal P_d)^{-1}
 =\xi^2(\xi^2+\mathcal L_d)^{-1}
  -i\xi\mathcal P_d(\xi^2+\mathcal L_d)^{-1}
 \qquad(\xi\in\mathbb R\setminus\{0\}).
\]
The first term is a normalised heat resolvent on the negative real axis. The second belongs to the preceding \(R\)-bounded family up to the unimodular factor \(-i\,\sgn(\xi)\). Since \(\mathcal P_d\) is sectorial of angle zero by {\rm(ii)}, this uniform imaginary-axis bound gives \ref{thm:main2} for the Poisson family by \cite[Proposition~10.3.3]{HNVW2}. Restricting to \(\lambda=-1/s\) gives \ref{thm:main2}\(\Rightarrow\)\ref{thm:main3}.

The chain \eqref{eq:Stein-six-bounds} and the bounds for \(C_{\mathrm{Poi}}\) identify the negative-axis \(R\)-bound of each family with \(\mathfrak s_{p,X}\) in \([0,\infty]\), proving \eqref{eq:resolvent-controls-Stein}. If \ref{thm:main3} holds for the chosen family, this common value is finite, so \ref{thm:main1} holds.
\end{proof}

\begin{corollary}[A one-dimensional resolvent characterisation]
\label{prop:onedim-Stein}
For every complex Banach space \(X\) and every \(1<p<\infty\),
\begin{equation*}
 \mathfrak s_{p,X}=C_{\mathrm{Poi},1}=C_{\Delta,1}
 \qquad\text{in }[0,\infty].
\end{equation*}
The \(R\)-bound of the Poisson resolvents on \(L^p(\mathbb T;X)\) also equals \(\mathfrak s_{p,X}\).
\end{corollary}

\begin{proof}
The first equality follows from \eqref{eq:Stein-six-bounds}. For the circle, Proposition~\ref{prop:onedim-circle-recovery} and \eqref{eq:onedim-real-to-circle} bound its \(R\)-bound from below by \(\mathfrak s_{p,X}\) and from above by \(C_{\mathrm{Poi},1}\).
\end{proof}

\section*{AI disclosure statement}

This manuscript originates from the observation that the proof of uniform \(R\)-sectoriality in \cite{Neerven-OU} uses the UMD assumption only through Stein's inequality and \(K\)-convexity. OpenAI's ChatGPT subsequently pointed out the known fact that Stein's inequality implies \(K\)-convexity. From here, the manuscript developed through extensive conversations with ChatGPT, which included exploratory mathematical discussions to find proofs, checking arguments, and polishing. The figure was generated with ChatGPT as well. In a later stage, Anthropic's Claude was used to review the manuscript and to suggest further improvements. The author has checked all arguments and takes full responsibility for the mathematical content and the final version of the paper.


\begin{thebibliography}{99}

\bibitem{BetancorEtAl}
J.~J. Betancor, A.~J. Castro, J.~Curbelo and L.~Rodr\'iguez-Mesa, \emph{Characterization of UMD Banach spaces by imaginary powers of Hermite and Laguerre operators}, Complex Anal. Oper. Theory \textbf{7} (2013), no.~4, 1019--1048. \href{https://doi.org/10.1007/s11785-011-0203-9} {doi:10.1007/s11785-011-0203-9}.
\bibitem{Bourgain83}

J.~Bourgain, \emph{Some remarks on Banach spaces in which martingale difference sequences are unconditional}, Ark. Mat. \textbf{21} (1983), no.~1, 163--168. \href{https://doi.org/10.1007/BF02384306}{doi:10.1007/BF02384306}.

\bibitem{Bourgain86}
J.~Bourgain, \emph{Vector-valued singular integrals and the \(H^1\)--BMO duality}, in: Probability Theory and Harmonic Analysis (Cleveland, Ohio, 1983), Monogr. Textbooks Pure Appl. Math. \textbf{98}, Dekker, New York, 1986, 1--19.

\bibitem{Eskenazis}
A. Eskenazis, \emph{On Pisier's inequality for UMD targets}, Canad. Math. Bull. \textbf{64} (2021), no. 2, 282--291. \href{https://doi.org/10.4153/S0008439520000442} {doi:10.4153/S0008439520000442}

\bibitem{Garling}
D.~J.~H. Garling, \emph{Random martingale transform inequalities}, in: Probability in Banach Spaces 6 (Sandbjerg, 1986), Progr.\ Probab.\ \textbf{20}, Birkh\"auser, Boston, 1990, 101--119. \href{https://doi.org/10.1007/978-1-4684-6781-9_6} {doi:10.1007/978-1-4684-6781-9\_6}.

\bibitem{Hytonen-LPS}
T.P.~Hyt\"onen, \emph{Littlewood--Paley--Stein theory for semigroups in UMD spaces},
Rev. Mat. Iberoam. \textbf{23} (2007), no.~3, 973--1009.
\href{https://doi.org/10.4171/RMI/521} {doi:10.4171/RMI/521}.

\bibitem{HytonenLiNaor}
T.P.~Hyt\"onen, S.~Li and A.~Naor, \emph{Quantitative affine approximation for UMD targets}, Discrete Anal. (2016), Paper No.~6, 48 pp. \href{https://doi.org/10.19086/da.614}{doi:10.19086/da.614}.

\bibitem{HNVW1}
T.P.~Hyt\"onen, J.M.A.M.~van Neerven, M.C.~Veraar and L.W.~Weis, \emph{Analysis in Banach Spaces, Volume I: Martingales and Littlewood--Paley Theory}, Ergebnisse der Mathematik und ihrer Grenzgebiete, vol. 63, Springer, 2016.
\href{https://doi.org/10.1007/978-3-319-48520-1}{doi:10.1007/978-3-319-48520-1}.

\bibitem{HNVW2}
T.P.~Hyt\"onen, J.M.A.M.~van Neerven, M.C.~Veraar and L.W.~Weis, \emph{Analysis in Banach Spaces, Volume II: Probabilistic Methods and Operator Theory}, Ergebnisse der Mathematik und ihrer Grenzgebiete, vol. 67, Springer, 2017. \href{https://doi.org/10.1007/978-3-319-69808-3} {doi:10.1007/978-3-319-69808-3}.

\bibitem{HNVW3}
T.P.~Hyt\"onen, J.M.A.M.~van Neerven, M.C.~Veraar and L.W.~Weis, \emph{Analysis in Banach Spaces, Volume III: Harmonic Analysis and Spectral Theory}, Ergebnisse der Mathematik und ihrer Grenzgebiete, vol. 76, Springer, 2023. \href{https://doi.org/10.1007/978-3-031-46598-7} {doi:10.1007/978-3-031-46598-7}.

\bibitem{KaltonLoristWeis}
N.~J. Kalton, E.~Lorist and L.W.~Weis, \emph{Euclidean Structures and Operator Theory in Banach Spaces}, Mem. Amer. Math. Soc. \textbf{288} (2023), no.~1433, vi+156 pp. \href{https://doi.org/10.1090/memo/1433} {doi:10.1090/memo/1433}.

\bibitem{KaltonWeis}
N.J.~Kalton and L.W.~Weis, \emph{The \(H^\infty\)-calculus and sums of closed operators}, Math. Ann. \textbf{321} (2001), no.~2, 319--345. \href{https://doi.org/10.1007/s002080100231}{doi:10.1007/s002080100231}.

\bibitem{Kato-Fractional-Note}
T.~Kato, \emph{Note on fractional powers of linear operators}, Proc. Japan Acad. \textbf{36} (1960), 94--96. \href{https://doi.org/10.3792/pja/1195524082}{doi:10.3792/pja/1195524082}.

\bibitem{Kato-Fractional}
T.~Kato, \emph{Fractional powers of dissipative operators}, J. Math. Soc. Japan \textbf{13} (1961), no.~3, 246--274. \href{https://doi.org/10.2969/jmsj/01330246}{doi:10.2969/jmsj/01330246}.

\bibitem{Komatsu-Fractional}
H.~Komatsu, \emph{Fractional powers of operators}, Pacific J. Math. \textbf{19} (1966), no.~2, 285--346. \href{https://doi.org/10.2140/pjm.1966.19.285}{doi:10.2140/pjm.1966.19.285}.

\bibitem{KwapienVeraarWeis}
S.~Kwapie\'n, M.C.~Veraar and L.W.~Weis, \emph{\(R\)-boundedness versus \(\gamma\)-boundedness}, Ark. Mat. \textbf{54} (2016), no.~1, 125--145. \href{https://doi.org/10.1007/s11512-015-0223-1} {doi:10.1007/s11512-015-0223-1}.

\bibitem{deLeeuw}
K.~de Leeuw, \emph{On \(L_p\) multipliers}, Ann. of Math. (2) \textbf{81} (1965), no.~2, 364--379.
\href{https://doi.org/10.2307/1970621}{doi:10.2307/1970621}.

\bibitem{Neerven-FA}
J.M.A.M.~van Neerven, \emph{Functional Analysis}, Cambridge Stud. Adv. Math. \textbf{201}, Cambridge Univ. Press, Cambridge, 2022. \href{https://doi.org/10.1017/9781009232487}{doi:10.1017/9781009232487}.

\bibitem{Neerven-OU}
J.M.A.M.~van Neerven, \emph{Dimension-free \(H^\infty\)-calculus of angle \(<\pi/2\) for UMD-valued Ornstein--Uhlenbeck operators}, preprint, 2026, \href{https://arxiv.org/abs/2609.04974} {arxiv.org/abs/2609.04974}.

\bibitem{Pisier}
G.~Pisier, \emph{Holomorphic semigroups and the geometry of Banach spaces}, Ann.\ of Math. \textbf{115} (1982), 375--392. \href{https://doi.org/10.2307/1971396} {doi:10.2307/1971396}.

\bibitem{Pisier-Martingales}
G.~Pisier, \emph{Martingales in Banach Spaces}, Cambridge Stud. Adv. Math. \textbf{155}, Cambridge Univ. Press, Cambridge, 2016.

\bibitem{Rota}
G.-C.~Rota, \emph{An \lq\lq Alternierende Verfahren'' for general positive operators}, Bull. Amer. Math. Soc. \textbf{68} (1962), no.~2, 95--102. \href{https://doi.org/10.1090/S0002-9904-1962-10737-X} {doi:10.1090/S0002-9904-1962-10737-X}.

\bibitem{Stein-LP}
E.M.~Stein, \emph{Topics in Harmonic Analysis Related to the Littlewood--Paley Theory}, Ann. of Math. Stud. \textbf{63}, Princeton Univ. Press, Princeton, NJ, 1970.


\end{thebibliography}
\end{document}